%%%%%%%%%%%%%%%%%%%%%%%%%%%%%%%%%%%%%%%%%%%%%%%%%%%%%%%%%%
%%%%%%%%%%%%%%%%%%%%%%%%%%%%%%%%%%%%%%%%%%%%%%%%%%%%%%%%%%
%%%    This is the AMS-LaTeX file:
%%
%%     Colli-Gilardi-Signori-Sprekels CGSS4
%%
%%       
%%%%%%%%%%%%%%%%%%%%%%%%%%%%%%%%%%%%%%%%%%%%%%%%%%%%%%%%%%
%%%%%%%%%%%%%%%%%%%%%%%%%%%%%%%%%%%%%%%%%%%%%%%%%%%%%%%%%%

\def\LaTeX{%
  \let\Begin\begin
  \let\End\end
  \let\salta\relax
  \let\finqui\relax
  \let\futuro\relax}

\def\UK{\def\our{our}\let\sz s}
\def\USA{\def\our{or}\let\sz z}

\UK 
%\USA

%%%%%%%%%%%%%%%%%%%%%%%%%%%%%%%%%

% scegliere fra \TeX e \LaTeX  e fra  \UK oppure \USA

%\TeX
\LaTeX

%\UK
\USA

%%%%%%%%%%%%%%%%%%%%%%%%%%%%%%%%%
%% page layout
%%%%%%%%%%%%%%%%%%%%%%%%%%%%%%%%%

\salta

\documentclass[twoside,12pt]{article}
\setlength{\textheight}{24cm}
\setlength{\textwidth}{16cm}
\setlength{\oddsidemargin}{2mm}
\setlength{\evensidemargin}{2mm}
\setlength{\topmargin}{-15mm}
\parskip2mm

%%%%%%%%%%%%%%%%%%%%%%%%%%%%%%%%%
%% packages
%%%%%%%%%%%%%%%%%%%%%%%%%%%%%%%%%

\usepackage{color}
\usepackage{amsmath}
\usepackage{amsthm}
\usepackage{amssymb,bbm}
\usepackage[mathcal]{euscript}

\usepackage{cite}
\usepackage{hyperref}
\usepackage[shortlabels]{enumitem}

\usepackage[ulem=normalem,draft]{changes}
%
%		COLORS FOR CORRECTIONS
%
% do the same, please (i.e., don't use the standard {\color{red} text} or similar): 
% just choose the color you prefer in \def\yourname

% EXAMPLE OF USE:  \fredi{I want this to become blue}
%
%IF YOU LATER WANT TO LET THE COLOR DISAPPEAR, ACTIVATE \def\fredi #1{{#1}} BELOW
 
%\definecolor{viola}{rgb}{0.3,0,0.7}
\definecolor{lilla}{rgb}{0.8,0,0.8}
\definecolor{blu}{rgb}{0,0,0.8}
\definecolor{rosso}{rgb}{0.85,0,0}

\def\juerg #1{{\color{black}#1}}
\def\an #1{{\color{black}#1}}
\def\gianni #1{{\color{black}#1}}
\def\pier #1{{\color{black}#1}}

\def\last #1{{\color{blue}#1}}
\def\pcol #1{{\color{red}#1}}

\def\pcol #1{#1}
\def\last #1{#1}
%\def\comm #1{#1 \ignorespaces}
%\def\gcomm #1{#1 \ignorespaces}

%%%%%%%%%%%%%%%%%%%%%%%%%%%%%%%%%
%% you may adjust the baseline
%%%%%%%%%%%%%%%%%%%%%%%%%%%%%%%%%

%\renewcommand{\baselinestretch}{0.975}

%%%%%%%%%%%%%%%%%%%%%%%%%%%%%%%%%
%% bibliographystyle
%%%%%%%%%%%%%%%%%%%%%%%%%%%%%%%%%

\bibliographystyle{plain}

%%%%%%%%%%%%%%%%%%%%%%%%%%%%%%%%%
%% environments
%%%%%%%%%%%%%%%%%%%%%%%%%%%%%%%%%

%
\newtheorem{theorem}{Theorem}[section]

\newtheorem{corollary}[theorem]{Corollary}

\finqui

\def\Beq{\Begin{equation}}
\def\Eeq{\End{equation}}
\def\Bsist{\Begin{eqnarray}}
\def\Esist{\End{eqnarray}}

\def\Bthm{\Begin{theorem}}
\def\Ethm{\End{theorem}}

\def\Brem{\Begin{remark}\rm}
\def\Erem{\End{remark}}

\def\Bcenter{\Begin{center}}
\def\Ecenter{\End{center}}
\let\non\nonumber

%%%%%%%%%%%%%%%%%%%%%%%%%%%%%%%%%
%% macros
%%%%%%%%%%%%%%%%%%%%%%%%%%%%%%%%%

% macro salvate

% sottosezioni non numerate

\def\step #1 \par{\medskip\noindent{\bf #1.}\quad}
\def\jstep #1: \par {\vspace{2mm}\noindent\underline{\sc #1 :}\par\nobreak\vspace{1mm}\noindent}

\def\aand{\quad\hbox{and}\quad}
\def\Lip{Lip\-schitz}
\def\Holder{H\"older}
\def\Frechet{Fr\'echet}

\def\lhs{left-hand side}
\def\rhs{right-hand side}

% versioni inglesi (UK) o americane (USA)

%\def\analyz {analy\sz}

% bold, cal, grass e mathop

\def\multibold #1{\def\arg{#1}%
  \ifx\arg\pto \let\next\relax
  \else
  \def\next{\expandafter
    \def\csname #1#1\endcsname{{\boldsymbol #1}}%
    \multibold}%
  \fi \next}

\def\pto{.}

\def\multical #1{\def\arg{#1}%
  \ifx\arg\pto \let\next\relax
  \else
  \def\next{\expandafter
    \def\csname cal#1\endcsname{{\cal #1}}%
    \multical}%
  \fi \next}

\def\multigrass #1{\def\arg{#1}%
  \ifx\arg\pto \let\next\relax
  \else
  \def\next{\expandafter
    \def\csname grass#1\endcsname{{\mathbb #1}}%
    \multigrass}%
  \fi \next}

% operatori

\def\multimathop #1 {\def\arg{#1}%
  \ifx\arg\pto \let\next\relax
  \else
  \def\next{\expandafter
    \def\csname #1\endcsname{\mathop{\rm #1}\nolimits}%
    \multimathop}%
  \fi \next}

\multibold
qweryuiopasdfghjklzxcvbnmQWERTYUIOPASDFGHJKLZXCVBNM.  % esclusa t per non cambiare \tt

\multical
QWERTYUIOPASDFGHJKLZXCVBNM.

\multigrass
QWERTYUIOPASDFGHJKLZXCVBNM.

\multimathop
diag dist div dom mean meas sign supp .

% accorpamenti di formule citate:
% uso  \accorpa {prima}{seconda}
%      \Accorpa\cs prima seconda (con il comodo blank anche dopo)
% NB: \Accorpa definisce \cs come l'accorpamento delle due citazioni
% e scrive sul file.log

\def\Accorpa #1#2 #3 {\gdef #1{\eqref{#2}--\eqref{#3}}%
  \wlog{}\wlog{\string #1 -> #2 - #3}\wlog{}}

% macro comode

\def\separa{\noalign{\allowbreak}}

\def\graffe #1{\mathopen\{#1\mathclose\}}
\def\<#1>{\mathopen\langle #1\mathclose\rangle}
\def\norma #1{\mathopen \| #1\mathclose \|}

\def\aeQ{\checkmmode{a.e.\ in~$Q$}}
\def\aet{\checkmmode{a.e.\ in~$(0,T)$}}
\def\aat{\checkmmode{for a.a.\ $t\in(0,T)$}}

\let\hat\widehat
\def\cpto{\,\cdot\,}

\def\Pi{\widehat\pi}

\def\iot {\int_0^t}
\def\ioT {\int_0^T}
\def\intQt{\int_{Q_t}}
\def\bintQt{\int_{Q^t}}
\def\intQ{\int_Q}
\def\iO{\int_\Omega}

\def\dt{\partial_t}
\def\dtt{\partial_t^2}
\def\dn{\partial_{\nn}}

\def\0{{\boldsymbol {0} }}

\let\emb\hookrightarrow

\def\checkmmode #1{\relax\ifmmode\hbox{#1}\else{#1}\fi}

% insiemi numerici

\let\erre\grassR
\let\enne\grassN

% spazi di funzioni a valori vettoriali su [0,T], [0,t], [0,s], [0,+\infty), [\delta,T]

% Come ricordare: in generale i simboli L H W  C da soli per gli spazi su (0,T)
% gli stessi raddoppiati per (0,+\infty)
% aggiunta di t o s al simbolo per (0,t) e (0,s)
% aggiunta di d al simbolo semplice o doppio per intervalli (\delta,T) e (\delta,+\infty)
% il simbolo C e i suoi derivati mettono le quadre anziche' le tonde

% Esempi   \L2V   \L\infty\Vp   \W{1,1}H   \C0H   \LL2V   \calC0\Vp   \Ld2V  \calCdH

\def\genspazio #1#2#3#4#5{#1^{#2}(#5,#4;#3)}
\def\spazio #1#2#3{\genspazio {#1}{#2}{#3}T0}
\def\spaziot #1#2#3{\genspazio {#1}{#2}{#3}t0}

\def\L {\spazio L}
\def\H {\spazio H}
\def\W {\spazio W}
\def\Lt {\spaziot L}

\def\C #1#2{C^{#1}([0,T];#2)}

% spazi di funzioni su \Omega, \Gamma, Q e \Sigma

\def\Lx #1{L^{#1}(\Omega)}
\def\Hx #1{H^{#1}(\Omega)}
\def\Wx #1{W^{#1}(\Omega)}

\def\LQ #1{L^{#1}(Q)}

\def\Ldue{\Lx 2}
\def\Linfty{\Lx\infty}

\def\Huno{\Hx 1}
\def\Hdue{\Hx 2}

\def\Liq{{L^\infty(Q)}}

% simboli in bold

% lettere greche

\let\theta\vartheta

\let\phi\varphi

\let\TeXchi\chi                         % new \chi, exactly on the baseline
\newbox\chibox
\setbox0 \hbox{\mathsurround0pt $\TeXchi$}
\setbox\chibox \hbox{\raise\dp0 \box 0 }
\def\chi{\copy\chibox}

% quadratino di fine dimostrazione

% abbreviazioni specifiche del lavoro

\def\vbar{\overline v}

\def\psibar{\overline\psi}
\def\phibar{\overline\phi}
\def\mubar{\overline\mu}

\def\phizbar{\overline\phiz}

\def\normaV #1{\norma{#1}_V}

\def\normaVp #1{\norma{#1}_*}

\def\phiz{\phi_0}
\def\wz{w_0}
\def\wu{w_1}

\def\Vp{{V^*}}

\def\CO{C_\Omega}
\def\kuno {{\kappa_1}}
\def\kdue {{\kappa_2}}

\def\J{{\cal J}}
\def\S{{\cal S}}
\def\Jred{{\J}_{\rm red}}
\def\Uad{{\cal U}_{\rm ad}}
\def\UR{{\cal U}_R}
\let\bconv\circledast

%%%%%%%%%%%%%%%%%%%%%%%%%%%%%%%%%%%%%%%%%%%%%%%%%%

\usepackage{amsmath}
\DeclareFontFamily{U}{mathc}{}
\DeclareFontShape{U}{mathc}{m}{it}%
{<->s*[1.03] mathc10}{}

\DeclareMathAlphabet{\mathscr}{U}{mathc}{m}{it}

\newcommand{\ov}[1]{\overline{#1}}

%%%%%%%%%%%%%%%%%%%%%%%%%%%%%%
\Begin{document}
%%%%%%%%%%%%%%%%%%%%%%%%%%%%%%%%%

%%%%%%%%%%%%%%%%%%%%%%%%%%%%%%%%%
%% front page
%%%%%%%%%%%%%%%%%%%%%%%%%%%%%%%%%

%
\title{
	Optimal temperature distribution for a nonisothermal Cahn--Hilliard system with source term
}
\author{}
\date{}
\maketitle
\Bcenter
\vskip-1.5cm
{\large\sc Pierluigi Colli$^{(1)}$}\\
{\normalsize e-mail: {\tt pierluigi.colli@unipv.it}}\\[0.25cm]
{\large\sc Gianni Gilardi $^{(1)}$}\\
{\normalsize e-mail: {\tt gianni.gilardi@unipv.it}}\\[0.25cm]
{\large\sc Andrea Signori$^{(2)}$}\\
{\normalsize e-mail: {\tt andrea.signori@polimi.it}}\\[0.25cm]
{\large\sc J\"urgen Sprekels$^{(3)}$}\\
{\normalsize e-mail: {\tt juergen.sprekels@wias-berlin.de}}\\[.5cm]
$^{(1)}$
{\small Dipartimento di Matematica ``F. Casorati'', Universit\`a di Pavia}\\
{\small and Research Associate at the IMATI -- C.N.R. Pavia}\\
{\small via Ferrata 5, I-27100 Pavia, Italy}\\[.3cm] 
$^{(2)}$
{\small Dipartimento di Matematica, Politecnico di Milano}\\
{\small via E. Bonardi 9, I-20133 Milano, Italy}
\\[.3cm] 
$^{(3)}$
{\small Department of Mathematics}\\
{\small Humboldt-Universit\"at zu Berlin}\\
{\small Unter den Linden 6, D-10099 Berlin, Germany}\\
{\small and}\\
{\small Weierstrass Institute for Applied Analysis and Stochastics}\\
{\small Mohrenstrasse 39, D-10117 Berlin, Germany}\\[10mm]
%\date{\an{======$>$[Last update \today]$<$======}}

\Ecenter
\Begin{abstract}
\noindent 
\juerg{In this note, we study the optimal control of a  nonisothermal phase field system of Cahn--Hilliard type that constitutes an 
extension of the classical Caginalp model for nonisothermal phase transitions with a conserved order parameter. \an{The system} couples a Cahn--Hilliard type equation with source term for the order parameter with the universal balance law of internal energy. In place of the standard 
Fourier form, the constitutive law of the heat flux is assumed in the form given by the theory developed by Green and Naghdi, which 
accounts for a possible thermal memory of the evolution. This has the consequence that the balance law of internal energy becomes
a second-order in time equation for the {\it thermal displacement} or {\it freezing index}, that is, a primitive with respect to time of
the temperature. Another particular feature of our system is the presence of the source term in the equation for the order parameter, which
entails additional mathematical difficulties because the mass conservation of the order parameter\an{, typical of the classic Cahn--Hilliard equation,} is no longer satisfied. 
In this paper, we \an{analyze} the case that the double-well potential driving the evolution of the phase transition is differentiable, either
(in the regular case) on the whole set of reals or (in the singular logarithmic case) on a finite open interval; nondifferentiable cases like the
double obstacle potential are excluded from the analysis.  
We prove the Fr\'echet differentiability of the control-to-state operator between suitable Banach spaces for both the regular and the 
logarithmic cases and establish the solvability of the corresponding adjoint systems in order to derive the associated first-order necessary
optimality conditions for the optimal control problem. Crucial for the whole analysis to work is the so-called ``strict separation property'', 
which states
that the order parameter attains its values in a compact subset of the interior of the effective domain of the nonlinearity. While this
separation property turns out to be  
generally valid for regular potentials in three dimensions of space, it can be shown for the logarithmic case only in two dimensions. 
}

\vskip3mm
\noindent {\bf Keywords:} Optimal control, nonisothermal Cahn--Hilliard equation, thermal memory, Cahn--Hilliard equation with source term, Cahn--Hilliard--Oono equation.

\vskip3mm
\noindent 
{\bf AMS (MOS) Subject Classification:} {
		35K55, % Nonlinear parabolic equations
        35K51, % Initial-boundary value problems for second-order parabolic systems		
		49J20, % Calculus of variations and optimal control; optimization;  Existence theories in calculus of variations and optimal control 
		49K20, % Calculus of variations and optimal control; optimization: Optimality conditions for problems involving partial differential equations 
%		49K40, % Calculus of variations and optimal control; optimization; Optimality conditions: Sensitivity, stability, well-posedness 
		49J50.  % Calculus of variations and optimal control; optimization;  Fréchet and Gateaux differentiability in optimization	
		}
\End{abstract}

\pagestyle{myheadings}
\newcommand\testopari{\sc Colli -- Gilardi -- Signori -- Sprekels}
\newcommand\testodispari{\sc Optimal Control of a nonisothermal Cahn--Hilliard type model}
\markboth{\testopari}{\testodispari}
%
%%%%%%%%%%%%%%%%%%%%%%%%%%%%%%%%%
%% very beginning
%%%%%%%%%%%%%%%%%%%%%%%%%%%%%%%%%

\section{Introduction}
\label{SEC:INTRO}
\setcounter{equation}{0}
\noindent 
\juerg{Let $\Omega\subset\erre^d$, $d\in\{2,3\}$, be some open, bounded, and connected set 
having a smooth boundary \gianni{$\Gamma:=\partial\Omega$} and the outward unit normal field~$\,\nn$. 
Denoting by $\dn$ the directional derivative in the direction of $\nn$, and putting, with a fixed final time $T>0$,
$$
  Q_t:=\Omega\times(0,t) \,\mbox{ and }\,\Sigma_t:=\gianni\Gamma\times(0,t)\,\mbox{ for $\,t\in(0,T] $, as well as \, $Q:=Q_T$ \,and
  \,$\Sigma:=\Sigma_T $},
$$ 
we study in this paper as \an{\it state system} the following initial-boundary value problem:}
\begin{alignat}{2}
  & \dt\phi - \Delta\mu + \gamma \phi
  = f  
  \quad && \text{in $\juerg{Q}$},
  \label{Iprima}
  \\
  & \mu
  = - \Delta\phi + F'(\phi) + a - b \dt w
  \quad && \text{in $Q$},
  \label{Iseconda}
  \\
  & \dtt w - \Delta(\kuno \dt w + \kdue w) + \lambda \dt\phi
  = u
  \quad && \text{in $Q$},
  \label{Iterza}
  \\
  & \dn\phi
  = \dn\mu
  = \dn(\kuno \dt w + \kdue w)
  = 0 
  \quad && \text{on $\juerg{\Sigma}$},
  \label{Ibc}
  \\
  & \phi(0) = \phiz, 
  \quad 
  w(0) = \wz,
  \quad 
  \dt w(0) = \wu 
  \quad && \text{in $\Omega.$}
  \label{Icauchy}
\end{alignat}
\Accorpa\Ipbl Iprima Icauchy
The {\it cost functional} under consideration \juerg{is given by}
\begin{align}
	\an{
	\J ((\phi, w), u)
	:= }{}
	&
%	\frac {\an{\alpha_1}}2 \ioT \norma{\phi - \phi_Q}^2
%	+ \frac {\an{\alpha_2}}2 \, \norma{\phi(T) - \phi_\Omega}^2
%	\non
%	\\ 
%	& {}
%	+ \frac {\an{\alpha_3}}2 \ioT \norma{w - w_Q}^2
%	+ \frac {\an{\alpha_4}}2 \, \norma{w(T) -  w_\Omega}^2
%	\non
%	\\ 
%	& {}
%	+ \frac {\an{\alpha_5}}2 \ioT \norma{\dt w -  w'_Q}^2
%	+ \frac {\an{\alpha_6}}2 \, \norma{\dt w(T) -  w'_\Omega}^2
%	+ \frac{\nu}{2} \ioT \norma{u}^2,
	\ \an{\frac {\an{\alpha_1}}2 \ioT\!\!\!\iO |{\phi - \phi_Q}|^2
	+ \frac {\an{\alpha_2}}2 \, \iO|{\phi(T) - \phi_\Omega}|^2}
	\non
	\\ 
	& \an{{}
	+ \frac {\an{\alpha_3}}2 \ioT\!\!\! \iO|{w - w_Q}|^2
	+ \frac {\an{\alpha_4}}2 \, \iO|{w(T) -  w_\Omega}|^2}
	\non
	\\ 
	& \an{{}
	+ \frac {\an{\alpha_5}}2 \ioT\!\!\!\iO|{\dt w -  w'_Q}|^2
	+ \frac {\an{\alpha_6}}2 \, \iO|{\dt w(T) -  w'_\Omega}|^2
	+ \frac{\nu}{2} \ioT\!\!\!\iO|{u}|^2,}
	\label{cost}
\end{align}
\juerg{with nonnegative constants $\an{\alpha_i}$, $1\le i\le 6$, which are not all zero, and where $\phi_\Omega, w_\Omega, w'_\Omega\in\Ldue$ and
$\phi_Q,w_Q,w'_Q\in L^2(Q)$ denote given target functions.}

For the control variable $u$, we choose as control space 
\gianni{%
\Beq
  \calU := \LQ\infty,
  \label{defU}
\Eeq
}%
and the \juerg{related set of admissible controls
is given} by
\begin{align}
	\label{Uad}
	\Uad : = \big\{ u \in {\cal U}: u_{\rm min} \leq u \leq u_{\rm max} \quad \aeQ \big\},
\end{align}
\juerg{where $u_{\rm min},u_{\rm max}\in\Liq$ satisfy $u_{\rm min}\le u_{\rm max}$ almost everywhere in $Q$.}

In summary, the control problem under investigation can be reformulated as follows:
\begin{align*}
\hspace*{-10mm}	{\boldsymbol {\rm (P)}}
	\quad 
	 \min_{u \in \Uad}\J ((\phi, w), u) \quad \text{subject to the constraint that $(\phi,\mu, w) $ solves \Ipbl.}
\end{align*}

\juerg{
The state system \eqref{Iprima}--\eqref{Icauchy} is a formal extension of the nonisothermal Cahn--Hilliard system introduced by Caginalp in \cite{Cag2}
to model the phenomenon of nonisothermal phase segregation in binary mixtures (see also \cite{Cag1, Cag3} and the derivation in \cite[Ex.~4.4.2, 
(4.44), (4.46)]{BS}); it corresponds to the Allen--Cahn counterpart analyzed 
recently in \cite{CSS3}. The unknowns in the state system have the following physical meaning: $\phi$ is a normalized difference between 
the volume fractions
of pure phases in the binary mixture (the dimensionless {\em order parameter} of the phase transformation, which should attain its values in 
the interval $[-1,1]$), 
$\mu$ is the associated 
{\em chemical potential}, and $\,w\,$ is the so-called {\em thermal displacement} (or 
{\em freezing index}), which is directly connected to the temperature $\theta$ (which in the case of the Caginalp model is actually 
a temperature difference) through the relation
\begin{align}
	\label{thermal_disp}	
	w (\cdot , t)  = w_0 + \iot \theta(\cdot, s) \,ds, \quad t \in[0,T].
\end{align}
Moreover, $\kuno$ and $\kdue$ in \eqref{Iterza} stand for prescribed positive coefficients related to the heat flux, which 
is here assumed in the Green--Naghdi form
(see~\cite{GN91,GN92,GN93,PG09})
\begin{equation}\label{flux}
\mathbf q=-\kuno \nabla (\dt w )- \kdue \nabla w \quad\mbox{where $\kappa_i>0$, $i=1,2$,}
\end{equation} 
which accounts for a possible previous
thermal history of the phenomenon. Moreover,
$\gamma$ is a positive physical constant related to the intensity of the mass absorption/production of the source,
where the source term in \eqref{Iprima} is $S:=f- \gamma \varphi$. This term reflects the fact that the system may not be isolated and the loss or production of mass is possible, which happens, e.g., in numerous liquid-liquid phase segregation problems that arise in cell biology~\cite{Bio}
and in tumor growth models~\cite{GLSS}. Notice that the presence of the source term entails that the property of mass conservation of the order
parameter is no longer valid; in fact, from \eqref{Iprima} it directly follows that the mass balance has the form
\begin{equation}
\label{masscons} 
\frac d{dt} \, \Bigl(\frac 1 {|\Omega|}\iO \phi(t) \Bigr)
	= 
	\frac 1 {|\Omega|} \iO {S(t)}, \quad\mbox{for a.e. $t\in(0,T)$},
\end{equation}
where $\,|\Omega|\,$ denotes the volume of $\Omega$.} \pier{To this concern, we would like 
to quote the paper \cite{CGRS}, where \last{a comparable} Cahn--Hilliard system without 
mass conservation was examined from the \last{optimal control viewpoint. 
Moreover, we refer \pcol{to \cite{CGS, CL, CRW, MTT, Per, Scar, SW}, 
where similar systems have been analyzed. For optimal control problems 
involving sparsity effects, let us mention \cite{CSS4, GaLaSi, SigTime, ST}.}}}
\an{Also, let us incidentally point out that the differential structure of equation 
\eqref{Iterza}, with respect to $w$, is sometimes also referred to as the {\it strongly 
damped wave equation}, see, e.g., \cite{PS} and the references therein.%
}

\juerg{%
In addition to the quantities already introduced,
$\lambda$ stands for the latent heat of the phase transformation, $a, b$ are physical constants, and the control variable
$\,u\,$ is  a distributed heat source/sink. Besides,  $\phi_0,w_0,$ and $w_1$ indicate some given initial values.  
Finally, the function $F$ is assumed to have a double-well shape.} Prototypical choices for the double-well shaped 
nonlinearity $F$ are the regular and singular {\it logarithmic potential} and its common (nonsingular) polynomial approximation, the {\it regular potential}.
In the order, they are defined as
\begin{align}
  & 
  \label{Flog}
  F_{log}(r) := 
  \left\{
    \begin{array}{ll}
      (1+r)\ln (1+r)+(1-r)\ln (1-r) - c_1 r^2 
      & \quad \hbox{if $|r|\leq1$},
      \\
      +\infty
      & \quad \hbox{otherwise},
	\end{array}
  \right.
  \\
  &   
	F_{reg}(r) := \frac 14 \, (r^2-1)^2 \,,
  \quad r \in \erre ,
  \label{Freg}
\end{align}
with the convention that $0\ln(0):=\lim_{r\searrow0}r\ln(r)=0$ and  $c_1>1$ so that $F_{log}$ is nonconvex.
Another important example is the nonregular and singular {\em double obstacle potential}, given~by
\Beq
  F_{2obs}(r) := - c_2 r^2 
  \quad \hbox{if $|r|\leq1$}
  \aand
  F_{2obs}(r) := +\infty
  \quad \hbox{if $|r|>1$},
  \label{F2ob}
\Eeq
with $c_2>0$. \juerg{However, the double obstacle case is not included in the subsequent analysis,} \pier{although we expect that,
\an{with similar techniques as those employed in \cite{CGSS2},} it is possible to extend the analysis also to this kind of nonregular potentials.}

\juerg{The state system \eqref{Iprima}--\eqref{Icauchy} was recently in \cite{CGSS3} analyzed concerning well-posedness and
regularity (see the results cited below in Section 2), where also the double obstacle case was included. Here, we concentrate on the
optimal control problem. While the existence of optimal controls is not too difficult to show, the derivation of 
first-order necessary optimality conditions is a much more challenging task, since it makes the derivation of differentiability 
properties of the associated control-to-state operator necessary. This, however, requires that the order parameter $\,\phi\,$ satisfies
the so-called \emph{strict separation property}, which means that $\,\phi\,$ attains its values in a compact subset of the
interior of the effective domain of the derivative $\,F'\,$ of $\,F$. While for regular potentials \an{this} condition turns out to be 
generally satisfied, \pier{it} cannot be guaranteed for singular potentials. In fact, following the ideas of the recent paper
\cite{CGRS4} on the isothermal case, one is just able to ensure the validity of the strict separation property for the logarithmic
potential $F_{\rm log}$ in the two-dimensional case $d=2$. Correspondingly, the analysis leading to first-order necessary optimality
conditions will be restricted to either the regular case for $d\le 3$ or the logarithmic case in two dimensions of space. In this sense, our
results apply to similar cases as those studied 	in \cite{CGRS4} in the isothermal situation. Observe, however, that the control problem considered
in \cite{CGRS4} differs considerably from that studied in this paper: indeed, in \cite{CGRS4} the control $\,u\,$ occurs in the 
order parameter equation resembling \eqref{Iprima}, while in our case it appears in the energy balance \eqref{Iterza}; for this reason, the set of 
admissible controls $\Uad$ had to be assumed in \cite{CGRS4} as a subset of the space $H^1(0,T;\Ldue)\cap\Liq$, which is cumbersome from
the viewpoint of optimal control, instead of the much
better space $\Liq$ used here. }

The plan of the paper is as follows.
The next section is devoted to collect previous results concerning the well-posedness of the state system \Ipbl.
Then, under suitable conditions, we provide some stronger \an{analytic} results in terms of regularity and stability properties of the state 
system with respect to the control variable $u$ appearing in \eqref{Iterza}.
The proof of \juerg{these} new results are addressed in Section \ref{SEC:REGCD}. Then, using \juerg{these results}, we analyze in 
Section \ref{SEC:CONTROL} the optimal control problem {\bf (P)}.

\section{Notation, assumptions and analytic results}
\label{SEC:RESULTS}
\setcounter{equation}{0}

First, let us set some notation and general assumptions.
\juerg{For any} Banach space $X$, we employ the notation $\norma{\cdot}_X$, $X^*$, and $\< \cdot , \cdot >_X$, to indicate 
the corresponding norm, its dual space, and  the related duality pairing between $X^*$ and~$X$.
For two Banach spaces $X$ and $Y$ \gianni{continuously embedded in some topological vector space~\an{$Z$}}, we introduce the linear space
$X\cap Y$, which becomes a Banach space when equipped with its natural norm 
\gianni{$\norma v_{X\cap Y}:=\norma v_X+\norma v_Y\,$}, for $v\in X\cap Y$. 

A special notation is used for the standard Lebesgue and Sobolev spaces defined on $\Omega$.
\juerg{For every $1 \leq p \leq \infty$ and $k \geq 0$, they} are denoted by $L^p(\Omega)$ and $W^{k,p}(\Omega)$, 
with the associated norms $\norma{\cpto}_{L^p(\Omega)}=\norma{\cpto}_{p}$ and $\norma{\cpto}_{W^{k,p}(\Omega)}$, respectively. 
\gianni{If $p=2$, they become Hilbert spaces, and we employ the standard convention $H^k(\Omega):= W^{k,2}(\Omega)$.} 
\juerg{For convenience}, we also~set
\begin{align*}
  & H := \Ldue , \quad  
  V := \Huno,   \quad
  W := \{v\in\Hdue: \ \dn v=0 \,\mbox{ on $\,\Gamma$}\}.
%  \label{def:HVW}
\end{align*}
\gianni{For simplicity, we use the symbol $\norma\cpto$ for the norm in $H$ and in any power of~it}.
\juerg{Observe that the embeddings} 
%\Beq
 $\, W \emb V \emb H \emb V^* \pier{{}\emb W^*{}}\,$
 % \non
%\Eeq
are dense and compact.
As usual, $H$ is identified with a subspace of $\Vp$ to have the Hilbert triplet $(V,H,\Vp)$ along with identity
\begin{align*}
	\text{$\<u,v>=(u,v)$ 
	\quad 
	for every $u\in H$ and $v\in V$,}
\end{align*}
\an{where we employ the special notation $\<\cdot,\cdot> := \<\cdot,\cdot>_V$. }

Next, for a generic element $v\in\Vp$,  we define its generalized mean value $\vbar$ by 
\Beq
  \vbar := \frac 1{|\Omega|} \, \< v , \boldsymbol{1} >,
  \label{defmean}
\Eeq
where $\boldsymbol{1}$ stands for the constant function that takes the value $1$ in $\Omega$.
It is clear that $\vbar$ reduces to the usual mean value if $v\in H$.
The same notation $\vbar$ is employed also if $v$ is a time-dependent function.

To conclude, for normed spaces $\,X\,$ and $\,v\in L^1(0,T;X)$, we define the convolution products
\begin{align}
	(\boldsymbol{1} * v)(t):=\int_0^t v(s)\,{\rm ds},
%	\label{conv}
	\quad 
	(\boldsymbol{1} \circledast v)(t):=\int_t^T v(s)\,{\rm ds} ,
	\qquad \hbox{$t \in[0,T]$.}
	\label{conv:back}	
\end{align}

For the \juerg{remainder of this paper, we make the following general assumptions}.
\begin{enumerate}[label={\bf (A\arabic{*})}, ref={\bf (A\arabic{*})}]
\item \label{ass:1:hpconst}
	The structural constants $\gamma$, $a$, $b$, $\kuno$, $\kdue$, and $\lambda$ are positive.
	
\item \label{ass:2:pot}
	The double-well potential $F$ can be written as $F= \hat \beta + \hat \pi$, where
	\begin{align*}
		\text{$\hat \beta:\erre\to [0,+\infty]$ is convex and lower semicontinuous with
	$\hat \beta(0)=0$}.
	\end{align*}
%	
%	$\hat \beta:\erre\to [0,+\infty]$ is convex and lower semicontinuous with
%	$\hat \beta(0)=0$, 
	This entails that $\beta := \partial\hat \beta$ is a maximal monotone graph with $\beta (0) \ni 0$.
	\juerg{Moreover, we assume that} 
\begin{align*}
	\text{$\hat \pi \in C^3(\erre)$, \juerg{where} $\pi:= \hat \pi ': \erre \to \erre$ is a \Lip\ continuous function.}
\end{align*}	
Besides, \juerg{denoting the effective domain of $\beta$ by $D(\beta)$, we assume that}  $D(\beta) = (r_-,r_+)\,$ with $\,-\infty\le r_-<0<r_+\le +\infty\,$ and \juerg{that} the restriction of $\hat \beta$ to $\,(r_-,r_+)\,$ belongs to $\,{C^3}(r_-,r_+)$.
There, $\beta$ reduces to \pier{the derivative of $\hat \beta$}, and we require that
\begin{align*}
	\text{$\lim_{r\searrow r_-} \beta(r)=-\infty$ \,and\, $\lim_{r\nearrow r_+}\beta(r)=+\infty$.}
\end{align*}
\pier{Please note that $F'$ in \eqref{Iseconda} has to be understood as $\beta + \pi$.}

%\item \label{ass:3:data}
%	For the data, we postulate that 
%\begin{align*}
%	\text{$f \in \LQ\infty, \quad   \phiz \in W \subset C^0({\ov Q}), \quad   \wz \in V$, \quad  and \quad $\wu \in H.$}
%\end{align*}

\item \label{ass:3:compat}
	\juerg{Let $f\in\Liq$}. We set $\rho := \frac {\norma f_\infty} \gamma$ and assume the compatibility \juerg{condition that all of the quantities} 
\begin{align*}
  & \juerg{\inf_{x\in\Omega}\phiz(x) , \ 
  \sup_{x\in\Omega}\phiz(x)}, \
  - \rho - (\phizbar)^- \,, \ \rho + (\phizbar)^+
	\quad \hbox{belong to the interior of $D(\beta)$},
\end{align*}
where $(\cdot)^+$  and $(\cdot)^-$ denote the positive and negative part functions, respectively. 
%
%\item \label{ass:5:Uad}
%	The functions $u_{\rm min} , u_{\rm max}$ belong to $ \cal U$ with $u_{\rm min} \leq u_{\rm max} \,\, \aeQ$.
%	
%\item \label{ass:6:pot}
%	$\an{\alpha_1},\dots,\an{\alpha_6},$ and $\nu$ are nonnegative constants, not all zero. 
%	
%\item \label{ass:7:target}
%	The target functions fulfill $\phi_Q ,w_Q,  w_Q' \in L^2(Q)$, $\phi_\Omega \in V, {w_\Omega}\in H$, and $w_\Omega' \in V$.

\end{enumerate}

The analysis of the above system \Ipbl\ has been the subject of investigation in \cite{CGSS3}.
There, weak and strong well-posedness has been addressed for general potentials and source terms. 
Since here we aim at solving the optimal control problem {\bf (P)}, we are forced to work under the framework of strong solutions. 
This, in particular, forces us to restrict the investigation \juerg{to differentiable potentials, more precisely, to either regular ones 
like \eqref{Freg} or, under the further restriction that $d=2$, to the logarithmic potential from \eqref{Flog}}. 
Since we are going to assume \ref{ass:1:hpconst}--\ref{ass:3:compat} in any case, we state the following results under \juerg{these} assumptions, even if some \juerg{of the} conditions may be relaxed (cf.~\cite{CGSS3}).

As a consequence of \cite[Thms.~2.2, 2.3, and 2.5]{CGSS3}, we have the following well-posedness result for the initial-boundary value problem \Ipbl.
\Bthm[Well-posedness of the state system]
\label{THM:WP}
Suppose that \ref{ass:1:hpconst}--\ref{ass:3:compat} hold true, and let the data of the system fulfill
\begin{align}
	\label{in:data:1}
	f & \in \H1\Vp,
	\quad 
	u \in \L2 H,
	\\
	\label{in:data:2}
	\phiz & \in \Hx3 \cap W , \quad
  \wz \in V, 
  \quad
  \wu \in V .
\end{align}
Then, there exists a unique solution $(\phi,\mu,w)$ to the system \Ipbl\ satisfying
\begin{align}
  & \phi \in \H1 V  \cap \L\infty {\Wx {2,6}}
  \gianni{\quad \hbox{with} \quad \beta(\phi) \in \L\infty{\Lx6}},
  \label{regphi}
  \\
  & \mu \in \L\infty V,
  \label{regmu}
  \\
  & w \in \H2 H \cap \pier{{} \C1V},
  \label{regw}
\end{align}
\Accorpa\Regsoluz regphi regw
\juerg{as well as} the estimate
\begin{align}
  & \norma\phi_{\H1 V  \cap \L\infty {\Wx {2,6}} }
  +\norma\mu_{\L\infty V}
  + \norma{\beta(\phi)}_{\L\infty{\Lx6}}
  \non
  \\
  &  \quad
  + \norma w_{\H2 H \cap \last{\C1 V} }
  \leq K_1\,,
  \label{wp:est}
\end{align}
with some constant $K_1>0$ that depends only on the structure of the system,
$\Omega$, $T$, and upper bounds for the norms of the data 
and the quantities related to assumptions \eqref{in:data:1}--\eqref{in:data:2}.
Besides, let $u_i \in \L2 H$, $i=1,2$, and let
$(\phi_i,\mu_i,w_i)$ be \juerg{the} corresponding solutions.
Then it holds that
\gianni{%
\begin{align}
	& \norma{\phi_1-\phi_2}_{\L\infty\Vp\cap\L2V}
  + \norma{w_1-w_2}_{\H1H\cap\L\infty V}
  \non
  \\
  & \quad 
   \leq K_2 \norma{\an{\boldsymbol{1}*}(u_1-u_2)}_{\L2H}\,,
%  \qquad
  \label{contdep}
\end{align}
}%
with some $K_2>0$ that depends only on the structure of the system,
$\Omega$, $T$, and an upper bound for the norms of $\beta(\phi_1)$ and $\beta(\phi_2)$ in~$\LQ1$.
\Ethm
Let us remark that, \pier{due to \eqref{regphi}}, the compact embedding $\Wx{2,6} \emb C^0(\ov \Omega)$, and 
\juerg{classical} compactness 
results (see, e.g., \cite[Sect.~8, Cor.~4]{Simon}), it follows that 
$\phi \in C^0(\ov Q)$.

\gianni{%
\Brem
\label{Variational}
The above well-posedness result \an{refers} to the natural variational form 
of the homogeneous Neumann problem for equation \eqref{Iprima},
due to the low regularity of $\mu$ specified in~\eqref{regmu}\an{.}
However, it is clear that, thanks to \eqref{regphi}\pier{, \ref{ass:3:compat}} and the elliptic regularity theory,
we also have that $\mu\in\L2W$, so that we actually can write \eqref{Iprima} in its strong~form. 
\pier{A similar consideration can be repeated for \last{the linear combination} $\kuno \dt w +\kdue w $
\last{in} \eqref{Iterza} \last{as you can find in the remark below}.} 
\Erem
}%

\pier{%
\Brem
\label{C1V}
We point \last{out} that the regularity $\C1V$ for \pcol{the variable} $w$ stated in \eqref{regw} does not  directly \last{follow} from \cite[Thms.~2.2, 2.3, 2.5]{CGSS3}, where just the regularity $\W{1,\infty}V$ was noticed. 
\last{This, however, can be deduced} with the help of \eqref{Iterza}, rewritten as the parabolic equation 
\begin{align}
&\label{pierinpiu}
%\frac1{\kuno} \dt (\kuno \dt w + \kdue w) - \Delta(\kuno \dt w + \kdue w) 
\last{\frac1{\kuno} \dt y - \Delta y = \pcol{f_w},}
%\\
%&\quad \hbox{with \rhs\ in }\, \L2H ,
\quad 
\last{\text{with $\, y:=\kuno \dt w + \kdue w \,$ \pcol{and  $\, f_w:= u -  \lambda \dt\phi +  \frac{\kdue}{\kuno} \dt w$},}}
\end{align}
\last{where, due to the previous results, it readily follows that $f_w \in \L2 H$.} \pcol{Note that $y$} satisfies \pcol{\eqref{pierinpiu} along with the} Neumann homogeneous boundary condition in \eqref{Ibc}\last{,} and the initial condition (cf.~\eqref{Icauchy})
$$ 
\last{y(0)={}}(\kuno \dt w + \kdue w) (0) = \kuno w_1 + \kdue w_0 \in V. 
$$
Then, by a straightforward application of the parabolic
regularity theory
(see, e.g., \cite{Baio1, LioMag}), it turns out that 
$$
\pcol{y ={}} \kuno \dt w + \kdue w \in \H1H\cap\C0V\cap\L2W \pcol{.}
$$
\pcol{At this point, it is not difficult to check} that $w\in \C1V$, whereas \pcol{we cannot infer the regularity} $w\in \H1W$ unless when $w(0) = w_0 \in W $.
\Erem
}%

As will be clear in the forthcoming Section  \ref{SEC:CONTROL}, 
the analytic framework encapsulated in Theorem \ref{THM:WP} 
does not suffice to rigorously prove the \Frechet\ differentiability of the solution operator 
associated with the system \Ipbl\ (cf.~Theorem \ref{THM:FRE}\  \pier{further on}) 
which is a key point to formulate the first-order necessary conditions for \gianni{optimality} addressed in Section~\ref{SUB:FOC}. 
For this reason, before entering the study of the optimal control \gianni{problem}~{\bf (P)}, 
we present some refined analytic results which are now possible 
by virtue of the more restricting condition we are
assuming on the potentials. 
In particular, a key regularity property to include singular and regular potentials in the analysis of the optimal control problem 
is the so-called {\it strict separation property} for the order parameter~$\phi$. 
\gianni{%
This \an{means} that
the values of $\phi$ are always confined in a compact subset of the interior of~$D(\beta)$.
Notice that, if $D(\beta)=\erre$, then the boundedness of $\phi$ that follows from the previous theorem
already guarantees this property.
For singular potentials, when $D(\beta)$ is an open interval, 
that means that the singularities of the potential at the end-points of $D(\beta)$ 
are not reached by $\phi$ at any time, meaning that the potential and its derivative actually are globally \Lip\ continuous functions.
The proof of the following result, sketched in Section~\ref{SEC:REGCD}, 
is derived with minor modifications arguing as done in \cite[Prop.~2.6]{CGRS4}.
It ensures both more regularity for the solution and the desired separation property 
in the important case of the logarithmic potential \eqref{Flog} in two dimensions.
}%

\begin{theorem}[Regularity and separation principle]
\label{THM:SEP}
Suppose that \ref{ass:1:hpconst}--\ref{ass:3:compat} hold, let $d=2$, and $F$ be the logarithmic potential defined in \eqref{Flog}.
Moreover, in addition to \eqref{in:data:1}--\eqref{in:data:2}, let \pier{$f$ and 
the auxiliary datum $\mu_0$} fulfill
\begin{align}
	\label{ass:ini:strong}
	f & \in \gianni{\H1 H},
	\quad
	\pier{\mu_0:= - \Delta \phiz + F'(\phiz) +a - bw_1 \in W}.
\end{align}
Then, the unique solution $(\phi,\mu,w)$ obtained from Theorem \ref{THM:WP} additionally enjoys the regularity properties
\begin{align} 
	\label{reg:strong}
	\dt \phi \in \L\infty H \cap \L2 W, 
	\quad 
	\mu \in L^\infty(Q),
	\quad 
	\beta(\phi) \in L^\infty(Q),
\end{align}
as well as 
\begin{align*}
	\norma{\dt \phi}_{\L\infty H \cap \L2 W}
	+ \norma{\mu}_{L^\infty(Q)}
	+ \norma{\beta(\phi)}_{L^\infty(Q)}
	\leq K_4,
\end{align*}
for \juerg{some} $K_4>0$ that depends only on the structure of the system, the initial data,
$\Omega$, and~$T$.
Furthermore, assume that 
\begin{align*}
	r_- < \min_{x \in \ov\Omega} \phiz (x) \leq \max_{x \in \ov\Omega} \phiz(x) < r_+.
\end{align*}
Then, the order parameter $\phi$ enjoys the strict separation property, that is, 
there exist \gianni{real numbers $r_*$ and~$r^*$} depending only on the structure of the system such that
\begin{align*}
		r_- < r_* \leq \phi(x,t)\leq r^* < r_+ \quad \text{ for $a.e. \,\,(x,t) \in Q$}.
\end{align*}
\end{theorem}
%
%Since we are assuming the singular $\beta$ to be regular in the interior of its domain and the separation property holds, we have the following.

\pier{%
\Brem
\label{Giorgini}
We point out that the regularity for $\mu $ in \eqref{reg:strong} is a consequence of the regularity $\mu \in \L\infty W$ and of the Sobolev embedding $W\last{\emb}L^\infty (\Omega)$, which holds up to the three-dimensional case. Notice also that a class of potentials slightly more general than the logarithmic one in \eqref{Flog} may be possibly  considered: for this aim we refer to \cite[Thm.~5.1]{GGM}\last{,} where a strict separation property has been derived in a suitable framework.
\Erem
}%

As a straightforward consequence of the above \gianni{results}, we have the following.

\begin{corollary}
\label{COR:SEP}
Suppose that \gianni{either $D(\beta)=\erre$ or that} the assumptions of Theorem \ref{THM:SEP} are fulfilled. 
Then, there exists a positive constant $K_5 $ just depending on 
\gianni{the structure and an upper bound for the norms of the data of the system} such~that  
\begin{align}
	\norma{\phi}_{L^\infty(Q)}
	+ \max_{i=\an{0,1,2,3}} \norma{F^{(i)}(\phi)}_{L^\infty(Q)}
	\leq K_5.
	\label{cor:sep}
\end{align}
\end{corollary}

With the above regularity improvement, we are now in a position to obtain 
a stronger continuous dependence estimate concerning the controls.
\begin{theorem}[Refined continuous dependence result]
\label{THM:CD:STRONG}
Suppose that \ref{ass:1:hpconst}--\ref{ass:3:compat} hold.
\gianni{Moreover}, assume that the \gianni{first and second derivatives of the potential $F$ are} \Lip\ continuous.
Consider $u_i \in \L2 H$, $i=1,2$, and let
\gianni{$(\phi_i,\mu_i,w_i)$, $i=1,2$, be the corresponding solutions}. 
Then, it holds that
\begin{align}
	& \norma{\phi_1-\phi_2}_{\H1 \Vp \cap \L\infty V\cap \L2 W}
  + \norma{\mu_1-\mu_2}_{\L2 V}
  \non
  \\
  & \quad 
   + \norma{w_1-w_2}_{\H2 \Vp \cap \W{1,\infty} {\last{V}} \cap \H1 {{ \last W}}}
  \leq K_6 \norma{u_1-u_2}_{\L2H}\pier{,}
 \label{cd:strong}
\end{align}
with \juerg{some} $K_6>0$ that depends only on the structure of the system,
$\Omega$, and $T$.
\end{theorem}
Notice that the above result \gianni{holds} for regular potentials 
both in dimensions two and three, as for \juerg{these} the \Lip\ continuity of $F'$ follows 
as a consequence of \gianni{Theorem~\ref{THM:WP}}.
On the other hand, the logarithmic potential can be considered just in dimension two 
as a consequence of the separation principle established by \gianni{Theorem~\ref{THM:SEP}}.
It is worth pointing out that the regularity improvement obtained in Theorem~\ref{THM:SEP} 
does not require more regularity of the control variable~$u$.
In particular, the strong well-posedness for the system is guaranteed 
for any control $u \in \L2 H$ \pier{(in which the control space ${\cal U}$ is embedded,} \gianni{see~\eqref{defU})}.

Let us conclude this section by collecting some useful tools that will be employed later on.
We often owe to the \gianni{Young, Poincar\'e and compactness inequalities}:
\begin{align}
  & ab \leq \delta a^2 + \frac 1{4\delta} \, b^2
  \quad \hbox{for every $a,b\in\erre$ and $\delta>0$},
  \label{young}
  \\[2mm]
  \separa
  & \normaV v
  \leq \CO \, \bigl( \norma{\nabla v} + |\vbar| \bigr)
  \quad \hbox{for every $v\in V$},
  \label{poincare}
  \\[2mm]
  \separa
  & \gianni{\norma v
  \leq \delta \, \norma{\nabla v} + C_{\Omega,\delta} \, \normaVp v
  \quad \hbox{for every $v\in V$ and $\delta>0$},
  \label{compact}}
\end{align}
\gianni{where $\CO$ depends only on~$\Omega$, $C_{\Omega,\delta}$ depends on~$\delta$, in addition,
and $\normaVp\cpto$ is the norm in $\Vp$ \juerg{to be introduced} below (see \eqref{normaVp})}.

\gianni{Next}, we recall an important  tool \juerg{which is} commonly used when working with problems connected to the Cahn--Hilliard equation.
Consider the weak formulation of the Poisson equation $-\Delta z=\psi$
\gianni{with homogeneous Neumann boundary conditions}. 
Namely, for a given $\psi\in\Vp$ (and not necessarily in~$H$), we consider the problem:
\begin{align}
  \hbox{\pier{find}} \quad 
	z \in V
  \quad \hbox{such that} \quad
  \iO \nabla z \cdot \nabla v
  = \< \psi , v >
  \quad \hbox{for every $v\in V$}.
  \label{neumann}
\end{align}
Since $\Omega$ is connected and regular, it is \juerg{well known} that the above problem admits a unique solution $z$ if and only if $\psi$ has zero mean value.
Hence, we can introduce the associated solution operator $\cal N$, which turns out to be \juerg{an} isomorphism between the following spaces, as
\begin{align}
  & \calN: \dom(\calN):=\graffe{\psi\in\Vp:\ \psibar=0} \to \graffe{z\in V:\ \ov z=0},
  \quad 
  {\cal N}: \psi \mapsto z,
  \label{defN}
\end{align}
where $z$  is the unique solution to \eqref{neumann} satisfying $\ov z=0$.
Moreover, it follows that the formula
\begin{align}
  \normaVp\psi^2 := \norma{\nabla\calN(\psi-\psibar)}^2 + |\psibar|^2
  \quad \hbox{for every $\psi\in\Vp$}
  \label{normaVp}
\end{align}
defines a Hilbert norm in $\Vp$ that is equivalent to the standard dual norm of $\Vp$.
From the above properties, one can obtain the following identities:
\begin{align}
  & \iO \nabla\calN\psi \cdot \nabla v
  = \< \psi , v >
  \quad  \hbox{for every $\psi\in\dom(\calN)$, $v\in V$},
  \label{dadefN}
  \\
  & \< \psi , \calN\zeta >
  = \< \zeta , \calN\psi >
  \quad \hbox{for every $\psi,\zeta\in\dom(\calN)$},
  \label{simmN}
  \\
  & \< \psi , \calN\psi > 
  = \iO |\nabla\calN\psi|^2
  = \normaVp\psi^2
  \quad \hbox{for every $\psi\in\dom(\calN)$},
  \label{danormaVp}
\end{align}
as well as 
\begin{align}
  \iot \< \dt v(s) , \calN v(s) > \, ds
  = \iot \< v(s) , \calN(\dt v(s)) > \, ds
  = \frac 12 \, \normaVp{v(t)}^2
  - \frac 12 \, \normaVp{v(0)}^2\,,
  \label{propN} 
\end{align}
which holds for every $t\in[0,T]$ and every \pier{$v\in\H1{\dom(\calN)}$}.
\Accorpa\PropN defN propN

Finally, without further reference later on, we are going to employ the following convention: the capital-case symbol $C$ is used to denote every constant
that depends only on the structural data of the problem such as
$\Omega$, $T$, $a$, $b$, $\kuno$, $\kdue$, $\gamma$, $\lambda$, the shape of the 
nonlinearities, and the norms of the involved functions. Therefore, its meaning may vary from line to line and even within the same line.
In addition, when a positive constant $\delta$ enters the computation, the related symbol $C_\delta$, in place of a general $C$, denotes constants
that \juerg{depend} on $\delta$, in addition.

\section{Regularity and continuous dependence results}
\label{SEC:REGCD}
\setcounter{equation}{0}

This section is devoted to the \juerg{proofs} of Theorem \ref{THM:SEP} and Theorem \ref{THM:CD:STRONG}. 
The first result is propedeutic to the second one which will play a key role in proving that the solution operator associated with the system enjoys some differentiability properties.

\gianni{%
\begin{proof}[Proof of Theorem \ref{THM:SEP}]
We can follow exactly the same argument \juerg{as that} used in \cite[Sect.~5.2]{CGRS4}
to prove the analogous result \cite[Prop.~2.6]{CGRS4}.
However, although we should perform the estimates in a rigorous way 
on a suitable discrete scheme designed on a proper approximating problem
as done in the quoted paper,
we proceed formally, for simplicity, by directly acting on problem \Ipbl,
and point out the few differences arising from the presence of the additional variable~$w$.
We differentiate both \eqref{Iprima} and \eqref{Iseconda} with respect to time
and test the resulting inequalities by $\dt\phi$ and $\Delta\dt\phi$, respectively.
If we sum up and integrate by parts and over~$(0,t)$, then a cancellation occurs, and we obtain~that
\begin{align}
  & \frac 12 \iO |\dt\phi(t)|^2
  + \gamma \intQt |\dt\phi|^2
  + \intQt |\Delta\dt\phi|^2
  \non
  \\
  & \quad= \frac 12 \iO |\dt\phi(0)|^2
  + \intQt \dt \juerg{f} \, \dt\phi
  + \intQt (\beta'+\pi')(\phi) \, \dt\phi \, \Delta\dt\phi
  - b \intQt \dtt w \, \Delta\dt\phi  \,.
  \non
\end{align}
This is the analogue of \cite[formula (5.16)]{CGRS4}
and essentially differs from it just for the presence of the last term.
%\an{In this direction, let us stress the last condition in \eqref{ass:ini:strong} 
%that we need to assume for the corresponding initial datum.}
However, this \juerg{term} can be easily dealt with by \juerg{using Young's} inequality
%\an{, the fact that $\dt f \in \L2 H$,} 
and the regularity of $w$ ensured by~\eqref{regw}.
\juerg{Indeed, we have that}
\Beq
  - b \intQt \dtt w \, \Delta\dt\phi
  \leq \frac 14 \intQt |\Delta\dt\phi|^2 
  + C \intQt |\dtt w|^2
  \leq \frac 14 \intQt |\Delta\dt\phi|^2 
  + C \,.
  \non
\Eeq
As the other terms can be treated as in the quoted paper,
we arrive at the analogue of \cite[formula (5.17)]{CGRS4},~i.e.,
\begin{align}
  & \norma{\dt\phi}_{\L\infty H}^2
  + \norma{\Delta\dt\phi}_{\L2H}^2
  \non
  \\[2mm]
  & \quad\leq C \, \bigl(
    \norma{\dt\phi(0)}^2
    + \norma{\beta'(\phi)}_{\L2{\Lx3}}^2 + 1
  \bigr)
  \, e^{C \, \norma{\beta'(\phi)}_{\L4{\Lx3}}^4} \,.
  \non
\end{align}
At this point, the new variable $w$ just enters the computation of~$\dt\phi(0)$. 
By still proceeding formally,
\pier{we recover the initial value for $\mu (0) = \mu_0$ from \eqref{Iseconda} at the time $t=0$, then, using the regularity of $\mu_0$ (and $f$)  stated in \eqref{ass:ini:strong}, 
we find out that  $$\dt\phi(0)  = f(0) + \Delta \mu_0 - \gamma \phiz \in H$$ from \eqref{Iprima}, also written for $t=0$. Then, we obtain that}
%and $\Delta\dt\phi(0)$, respectively, and sum up.
%Since the terms involving $\mu(0)$ cancel each other, 
%with the help of some integrations by parts and our assumptions on the data, 
%we obtain~that
\begin{align}
   \pier{\norma{\dt\phi(0)}^2 \leq \norma{f(0) + \Delta \mu_0 - \gamma \phiz}^2 \leq C.}
%   &= \iO |\dt\phi(0)|^2
%  \non
%  \\
%  & {} = \iO \bigl( f(0) - \gamma \phi(0) \bigr) \, \dt\phi(0)
%  + \iO \nabla\phi(0) \cdot \nabla\Delta\dt\phi(0)
%  \non
%  \\
%  & \quad {}
%  + \iO \bigl( F'(\phi(0)) + a - b \dt w(0) \bigr) \Delta\dt\phi(0)
%  \non
%  \\
%  & {} = \iO \bigl\{
%      f(0) - \gamma \phiz
%      - \Delta^2 \phiz
%      + \Delta \bigl( F'(\phiz) + a - b \wu \bigr) 
%  \bigr\} \dt\phi(0)
%  \non
%  \\[1mm]
%  &{} \leq C \, \norma{\dt\phi(0)}\,, 
  \non
\end{align}
%whence the desired estimate of $\norma{\dt\phi(0)}$ \juerg{follows}.
At this point, $w$ does not enter the argument any longer, 
and we can proceed and then conclude \an{exactly} as in~\cite{CGRS4}.
\end{proof}
}%

\begin{proof}[Proof of Theorem \ref{THM:CD:STRONG}]

To begin with, let us set the following notation for the \gianni{differences involved in the statement}:
\begin{align*}
	\phi:=\phi_1-\phi_2,
	\quad 
	\mu:= \mu_1-\mu_2,
	\quad 
	u:= u_1- u_2,
	\quad 
	w:= w_1 - w_2.
\end{align*}
Next, we write the system \an{solved by} the differences that, in its strong form, reads~as
\begin{alignat}{2}
  & \dt\phi - \Delta\mu + \gamma \phi
  = 0
  \quad && \text{in $Q$,}
  \label{Iprima:diff}
  \\
  & \mu
  = - \Delta\phi + (F'(\phi_1) - F'(\phi_2))  - b \dt w
  \quad && \text{in $Q$},
  \label{Iseconda:diff}
  \\
  & \dtt w - \Delta(\kuno \dt w + \kdue w) + \lambda \dt\phi
  = u
  \quad && \text{in $Q$},
  \label{Iterza:diff}
  \\
  & \dn\phi
  = \dn\mu
  = \dn(\kuno \dt w + \kdue w)
  = 0 
  \quad && \text{on $\Sigma$},
  \label{Ibc:diff}
  \\
  & \phi(0) = 
  w(0) = 
  \dt w(0) = 0 
  \quad && \text{in $\Omega.$}
  \label{Icauchy:diff}
\end{alignat}

\step First estimate

First, we recall that \gianni{$F'$} is now assumed to be \Lip\ continuous.
Then, testing \eqref{Iprima:diff} by~$\phi$,  \eqref{Iseconda:diff} by~$\mu$, and adding the resulting equations \gianni{lead} us~to
%\begin{align*}
%	& \frac 12 \iO |\phi(t)|^2
%	+ \gamma \int_{Q_t} |\phi|^2
%	+ \int_{Q_t} |\mu|^2
%	= \int_{Q_t} (F'(\phi_1)-F'(\phi_2)) \mu
%	- b \int_{Q_t} \dt w \mu
%	\\ & \quad 
%	\leq 
%	\frac 12 \int_{Q_t} |\mu|^2
%	+ C\int_{Q_t} |\phi|^2
%	+ C \int_{Q_t} |\dt w|^2.
%\end{align*}
\begin{align*}
	& \frac 12 \, \frac d {dt} \, \norma{\phi}^2
	+ \gamma \norma{\phi}^2
	+ \norma{\mu}^2
	= \iO (F'(\phi_1)-F'(\phi_2)) \mu
	- b \iO \dt w \mu
	\\ & \quad 
	\gianni{{} \leq 
	\frac 12 \, \norma{\mu}^2
	+ C \bigl( \norma{\phi}^2+\norma{\dt w}^2 \bigl).}
\end{align*}
Now, recalling the continuous dependence estimate already proved in Theorem \ref{THM:WP}, we infer, after integrating over time, that 
\begin{align}
	\label{cd:1}
	\norma{\phi_1-\phi_2}_{\L\infty H}
	+ \norma{\mu_1-\mu_2}_{\L2 H}
	\leq 
	C \norma{1 * (u_1-u_2)}_{\L2 H}.
\end{align}

\step Second estimate

\gianni{First, \an{let us} \juerg{establish} an auxiliary estimate.
Since $F'$ and $F''$ are \an{supposed} to be \Lip\ continuous
and \eqref{wp:est} ensures a uniform bound for~$\norma{\nabla\phi_2}_\infty$, 
we have\an{, almost everywhere in $(0,T)$,} that
\Bsist
  && \normaV{F'(\phi_1)-F'(\phi_2)}
  \leq \norma{F'(\phi_1)-F'(\phi_2)}
  + \norma{F''(\phi_1)\nabla\phi_1 - F''(\phi_2)\nabla\phi_2}
  \non
  \\ 
  && \quad \leq C \, \norma\phi
  + \norma{F''(\phi_1)\nabla\phi}
  + \norma{(F''(\phi_1)-F''(\phi_2))\nabla\phi_2}
  \non
  \\ 
  && \quad \leq C \, \norma\phi
  + C \, \norma{\nabla\phi}
  \leq C \, \normaV\phi \,.
  \non
\Esist
Next, we multiply \eqref{Iprima:diff} by $1 / |\Omega|$ to obtain that
\Beq
  \frac d{dt} \, \phibar(t) 
  + \gamma \, \phibar(t)= 0
  \quad \aat,
  \label{meanphi:diff}
\Eeq
which entails that $\phibar(t)=0$ for every $t\in[0,T]$ since $\phibar(0)=0$.
In particular, besides~$\phi$, even} $\dt\phi $ has zero mean value. 
Thus, we are allowed to test \eqref{Iprima:diff} by ${\cal N} (\dt \phi)$, 
\eqref{Iseconda:diff} by $-\dt\phi$, and \eqref{Iterza:diff} by $\frac b\lambda\dt w$, 
and add the resulting identities.
By also accounting for the \Lip\ continuity of~$F'$ and the Young inequality, we deduce~that\juerg{, a.e. in $(0,T)$,} 
\begin{align*}
	& \norma{\dt \phi}_*^2
	+ \frac \gamma 2 \, \frac d {dt} \, \norma{\phi}_*^2
	+ \frac 12 \, \frac d {dt} \, \norma{\nabla \phi}^2
	+ \gianni{\frac b{2\lambda}} \, \frac d {dt} \, \norma{\dt w}^2
	+ \frac{\kuno b}{2\lambda} \, \norma{\nabla (\dt w)}^2
	+ \frac {\kdue b} {2 \lambda} \, \frac d {dt} \, \norma{\nabla w}^2
	\\ &  \quad 
	{} = \iO  (F'(\phi_1) - F'(\phi_2)) \dt \phi
	+ \frac b\lambda \iO  u \, \dt w
	\\ & \quad 
	\gianni{{} \leq C \bigl( \normaV\phi \, \normaVp{\dt\phi} + \norma u \, \norma{\dt w} \bigr)}
	\gianni{{} \leq \frac 12 \, \normaVp{\dt\phi}^2
	+ C \bigl( \normaV\phi^2 + \norma u^2 + \norma{\dt w}^2 \bigr).}
\end{align*}
\gianni{Hence, integrating over time and using \eqref{contdep}}\juerg{, we may} conclude that 
\begin{align}
	\non
	& \norma{\phi_1-\phi_2}_{\H1 \Vp \cap \L\infty V}
	+ \norma{w_1-w_2}_{\W{1,\infty} H \cap \H1 V}
	\\
	& \quad {}
	\leq 
	C\norma{u_1-u_2}_{\L2 H}.
	\label{cd:2}
\end{align}

\step Third estimate

\gianni{%
By testing \eqref{Iprima:diff} by $\mu$, we have that
\Beq
  \iO \dt\phi \, \mu
  + \iO |\nabla\mu|^2
  + \gamma \iO \phi \mu 
  = 0 \,.
  \non
\Eeq
Now, we recall that $\phi$ and $\dt\phi$ have zero mean value.
Hence, by also accounting for the Poincar\'e inequality~\pier{\eqref{poincare}}, we deduce~that
\Beq
  \iO |\nabla\mu|^2
  = - \iO \dt\phi \, (\mu-\mubar)
  \juerg{\,-\,} \gamma \iO \phi (\mu-\mubar)
  \leq \frac 12 \iO |\nabla\mu|^2
  + C \bigl( \normaVp{\dt\phi}^2 + \norma\phi^2 \bigr).
  \non
\Eeq
Therefore,} \pier{thanks to \eqref{cd:1} and \eqref{cd:2}}, it readily follows that
\begin{align}
	\label{cd:3}
	\pier{\norma{\mu_1- \mu_2}_{\L2 V}}
	\leq 
	C \norma{u_1-u_2}_{\L2 H}.
\end{align}

\step  Fourth estimate

A simple comparison argument in \eqref{Iseconda:diff}, along with the above estimates  and elliptic regularity theory, entails that 
\begin{align}
	\label{cd:4}
	\norma{\phi_1- \phi_2}_{\L2 W}
	\leq 
	C \norma{u_1-u_2}_{\L2 H}.
\end{align}

\step  Fifth estimate

\gianni{We take an arbitrary $v\in\L2V$, multiply \eqref{Iterza:diff} by $v$,
and integrate over $Q$ and by parts.
By rearranging and estimating, we easily obtain~that
\Beq
  \intQ \dtt w \, v
  \leq C \bigl(
    \norma u_{\L2H}
    + \norma{\dt w}_{\L2V}
    + \norma w_{\L2V}
    + \norma{\dt\phi}_{\L2\Vp}
  \bigr) \norma v_{\L2V} \,.
  \non
\Eeq
On account of the previous estimates, we conclude that}
\begin{align}
	\label{cd:5}
	\norma{\dtt w_1 - \dtt w_2}_{\L2 \Vp}
	\leq 
	C \norma{u_1-u_2}_{\L2 H}.
\end{align}

\step \last{Sixth estimate}

\last{
Arguing as in Remark \ref{C1V}, we now \pcol{rewrite} \eqref{Iterza:diff} as a parabolic equation in the auxiliary variable $y:= \kuno \dt w + \kdue w + \kuno \lambda \phi$ obtaining that
\begin{align*}
&
%\frac1{\kuno} \dt (\kuno \dt w + \kdue w) - \Delta(\kuno \dt w + \kdue w) 
\last{\frac1{\kuno} \dt y - \Delta y }
=  \pcol{u +  \displaystyle\frac{\kdue}{\kuno} \dt w - \kuno \lambda \Delta \phi.}
%\last{=:f_w.}
\end{align*}
%From the above estimates, we deduce that $\norma{f_w }_{\L2 H} \leq C \norma{u_1-u_2}_{\L2 H}$.
Besides, \pcol{let us underline that} $y$ satisfies homogeneous Neumann boundary \pcol{conditions and null initial conditions, as it can be realized from \eqref{Ibc:diff} and \eqref{Icauchy:diff}.
Then, using a well-known parabolic regularity result and the already found estimates \eqref{cd:2} and \eqref{cd:4}\last{,} it is straightforward to deduce that 
\begin{align*}
\norma{y}_{\H1H\cap\C0V\cap\L2W} \leq C \Bigl\| u +  \displaystyle\frac{\kdue}{\kuno} \dt w - \kuno \lambda \Delta \phi \Bigr\|_{\L2H}\, \leq C \norma u_{\L2H}\,.
\end{align*}
Thus, by solving the Cauchy problem for the ordinary differential equation $\kuno \dt w + \kdue w = y - \kuno \lambda \phi$ in terms of $w$, and recalling again \eqref{cd:2} and \eqref{cd:4}, we find out that} 
\begin{align}
	\label{cd:6}
	\norma{ w_1 - w_2}_{\pcol{\H2 \Vp \cap{}} \W {1,\infty} V \cap \H1 W}
	\leq 
	C \norma{u_1-u_2}_{\L2 H}.
\end{align}
}
This \gianni{completes} the proof, as collecting the above estimates \juerg{yields} \eqref{cd:strong}.
\end{proof}

\section{The optimal control problem}
\label{SEC:CONTROL}
\setcounter{equation}{0}

\juerg{In this section, we study} the optimal control problem introduced at the beginning\juerg{, which} we recall here for \juerg{the} reader's convenience:
\begin{align*}
	{\boldsymbol {(\rm P)}}
	\quad 
	\min_{u \in \Uad}\J ((\phi, w), u) \quad \text{subject to the constraint that $(\phi,\mu, w) $ solves \Ipbl,}
\end{align*}
\gianni{where the cost functional $\J$ is given by~\eqref{cost}.}

To begin with, let us fix some \juerg{notation} concerning the solution operator $\S$ associated with the system \Ipbl. 
As a consequence of the Theorems \ref{THM:WP}, \ref{THM:SEP}, and \ref{THM:CD:STRONG}, \juerg{the} {\it control-to-state operator} 
\begin{align*}
	\S \gianni{{} = (\S_1,\S_2,\S_3)} : L^2(Q) \pier{{}(\supset{\cal U}){}} \to {\cal Y}, \quad 
	\S: u \mapsto (\phi, \mu, w),
\end{align*}
\juerg{is well defined},
where $(\phi,\mu,w) \in \cal Y$ is the unique solution to the \juerg{state} system, and the Banach space $\cal Y$, referred to as the {\it state space}, is defined by the \gianni{regularity specified in \eqref{regphi}--\eqref{regw} and partially in \eqref{reg:strong}},
that~is,
\begin{align*}
	{\cal Y}&  := \big( \W{1,\infty} H  \cap \H1{W} \cap \L\infty {\Wx{2,6}}\big) \times \L\infty V
	 \\ & \quad 
	\times \big ( \H2 H \cap 
%	\W{1,\infty} V 
	\last{\C1 V}
	\big).
\end{align*}
Moreover, the \pier{continuous dependence} estimate provided by Theorem \ref{THM:CD:STRONG} 
allows us to infer that the solution operator is \Lip\ \juerg{continuous} in the sense that, 
\gianni{for any pair $(u_1,u_2)$ of controls}, it holds that
\begin{align*}
	& \norma{\S(u_1)- \S(u_2)}_{\cal X}
	\leq 
	K_6\norma{u_1-u_2}_{\L2H},
\end{align*}
where 
%with $\S_i$, $i=1,\dots,3$, we indicate its $i$-th component and 
$\cal X$ is the space defined by 
\begin{align}
	\non
	{\cal X} & := \big(\H1 \Vp \cap \L\infty V \cap \L2 W)\times \L2 V
	\\ & \quad \label{def:X}
	\times \big(\H2 \Vp \cap \pier{\W{1,\infty} V \cap \H1 W}).
\end{align}
Furthermore, we \juerg{introduce} the {\it reduced cost functional}\juerg{, given by}
\begin{align}\label{Jred}
	\Jred : L^2(Q) \subset{\cal U} \to \erre,
	\quad 
	\Jred: u \mapsto \J(\S_1( u),\S_3( u), u),
\end{align}
which allows us to reduce the optimization problem {\bf (P)} \juerg{to} the form
\begin{align*}
	\min_{u \in \Uad} \Jred (u).
\end{align*}

In \gianni{what} follows, we are working \juerg{in} the framework of Theorem \ref{THM:WP} 
(and \gianni{possibly} in the sense of Theorem \ref{THM:SEP}). 
For this reason, the following conditions will be in order: 
\begin{enumerate}
[label={\bf (C\arabic{*})}, ref={\bf (C\arabic{*})}]
\item \label{ass:4:data}
	The source $f$ \juerg{fulfills} \eqref{in:data:1}, and the initial data $\phiz, w_0,$ and \gianni{$w_1$} \juerg{satisfy}  \eqref{in:data:2}. 
	Moreover, if we consider the logarithmic potential and $d=2$, they additionally
	\gianni{fulfill} \eqref{ass:ini:strong}.

\item \label{ass:5:Uad}
	The functions $u_{\rm min} , u_{\rm max}$ belong to $ \cal U$ with $u_{\rm min} \leq u_{\rm max} \,\, \aeQ$.
	
\item \label{ass:6:pot}
	$\an{\alpha_1},\dots,\an{\alpha_6},$ and $\nu$ are nonnegative constants, not all zero. 
	
\item \label{ass:7:target}
	The target functions fulfill $\phi_Q ,w_Q,  w_Q' \in L^2(Q)$, $\phi_\Omega \in V, {w_\Omega}\in H$, and $w_\Omega' \in V$.
\end{enumerate}

\subsection{Existence of optimal controls}

The first result we are going to address concerns the existence of optimal controls.
\begin{theorem}[Existence of optimal controls] \label{THM:EX:OC}
We suppose that \juerg{the} \gianni{assumptions} \ref{ass:1:hpconst}--\ref{ass:3:compat} 
%\allowbreak
and \ref{ass:4:data}--\ref{ass:7:target} are fulfilled. 
Then, the optimal control problem {\bf (P)} admits a solution.
\end{theorem}

\begin{proof}[Proof of Theorem \ref{THM:EX:OC}]
As the proof is an immediate consequence of the direct method of the calculus of variations, we just briefly outline the crucial steps.
Consider a minimizing sequence $\{u_n\}_n \subset \Uad$ for the reduced  cost functional $\Jred$ defined by \eqref{Jred}.
Let us introduce also the sequence of the associated states $\{(\phi_n,\mu_n, w_n)\}_n$, where $(\phi_n,\mu_n, w_n) =\S(u_n)$ for every $n \in \enne$.
Namely, we have that
\begin{align*}
	\lim_{n \to \infty} \Jred(u_n) 
	= 	\lim_{n \to \infty} \J \big((\S_1(u_n), \S_3(u_n)), u_n \big) 
	= \inf_{u \in \Uad} \Jred(u)
	\geq 0.
\end{align*}
Thus, as $\Uad$ is bounded in $\cal U$, by standard compactness arguments, using also that $\Uad$ is closed and convex, we obtain a limit function $u^* \in \Uad$ and a nonrelabelled subsequence \juerg{such} that, as $n \to \infty$,
\begin{align*}
	u_n \to u^* \quad \text{weakly-star in $L^\infty(Q)$}.
\end{align*}
On the other hand, by the \pier{boundedness property \eqref{wp:est} stated in}
%continuity property of the solution operator which follows from 
Theorem \ref{THM:WP}, along with 
\juerg{standard} compactness results (see, e.g., \cite[Sect.~8, Cor.~4]{Simon}), we also \gianni{have} that
\begin{alignat*}{2}
	\phi_n & \to \phi^* \quad && \text{weakly-star in $\H1 V \cap \L\infty {\Wx{2,6}}$},
	\\ &  &&  \quad \text{and strongly in $C^0(\ov Q),$}
	\\
	\mu_n & \to \mu^* \quad && \text{weakly-star in $\L\infty V$},
	\\
	w_n & \to w^* \quad && \text{\gianni{weakly-star} in $\H2 H \cap \W{1,\infty} V$},
	\\ &  &&  \quad \text{and strongly in \gianni{$\C1H$},}
	\\
	\an{F'(\phi_n)} & \an{\to \xi^* \quad} && \an{\text{{weakly-star} in $\L\infty {\Lx6}$},}
\end{alignat*}
\an{for some limits $\phi^*,\mu^*, w^*$, and $\xi^*$.}
The first strong convergence follows from the compact embedding $\Wx {2,6} \emb C^0(\ov \Omega)$.
\an{Besides, \pier{as 
$$\pi(\phi_n) \to \pi(\phi^*) \quad \hbox{strongly in $C^0(\ov Q)$ \ \ and} 
\quad \beta (\phi_n) \to \xi^* - \pi (\phi^*) \quad \hbox{weakly in }L^1(Q), $$ 
by maximal monotonicity arguments}
%using standard arguments along with the above weak and strong convergences, 
it is not difficult to conclude that $\xi^* = F'(\phi^*)$.}
Then, using the above weak, weak-star and strong \gianni{convergence properties}, 
it is a standard matter to pass to the limit as $n$ \juerg{tends} to infinity in the variational formulation associated with system \Ipbl,
written for $\phi_n,\mu_n, w_n,$ and $ u_n$. This will \juerg{also prove} that $(\phi^*, \mu^*, w^*)$ is nothing but $\S(u^*)$.
Finally, the lower semicontinuity of norms entails that 
\begin{align*}
		\Jred (u^*) \leq \liminf_{n\to \infty}   \Jred (u_n ) = \lim_{n\to \infty}   \Jred (u_n ) = \inf_{ u\in \Uad} \Jred(u)\,,
\end{align*}
meaning that $u^* $ is a global minimizer for $\Jred$.
\end{proof}

\subsection{Differentiability \juerg{of} the solution operator}

In \juerg{the following,} we are going to prove some differentiability properties for the solution operator $\S$.
Since \juerg{these} have to be analyzed in open sets, let us \pier{take} an open ball in the $L^\infty$-topology \juerg{that}  contains 
\an{the set of admissible controls} $\Uad,$
namely, let $R>0$ \juerg{be chosen} such that
\begin{align*}
	  \juerg{ \UR := \{ u \in {\cal U} : \norma{u}_{\cal U} <R\} \supset  \Uad.}
\end{align*}
Now, we fix $u \in \UR$ and denote by $(\phi, \mu, w )= \S(u)$ the unique corresponding state. 
The\an{n, the} linearized system to \Ipbl\ at the fixed control $u$ is given, for any $h \in L^2(Q)$, as follows:
\begin{alignat}{2}
  & \dt\xi - \Delta\eta + \gamma \xi
  = 0
  \quad && \text{in $Q$},
  \label{lin:sys:1}
  \\
  & \eta
  = - \Delta\xi + F''(\phi)\xi - b \dt \zeta
  \quad && \text{in $Q$},
  \label{lin:sys:2}
  \\
  & \dtt \zeta - \Delta(\kuno \dt \zeta + \kdue \zeta) + \lambda \dt\xi
  = h
  \quad && \text{in $Q$},
  \label{lin:sys:3}
  \\
  & \dn\xi
  = \dn\eta
  = \dn(\kuno \dt \zeta + \kdue \zeta )
  = 0 
  \quad && \text{on $\Sigma$},
  \label{lin:sys:4}
  \\
  & \xi(0) =  \zeta(0) =  \dt\zeta (0) = 0
  \quad && \text{in $\Omega.$}
  \label{lin:sys:5}
\end{alignat}
\Accorpa\Lin {lin:sys:1} {lin:sys:5}
The \juerg{proof of the} well-posedness of the above system is very similar (and, in fact, easier, as the system is linear) to the proof of 
Theorem \ref{THM:WP}. \juerg{We have the following result.}

\begin{theorem}[Well-posedness of the linearized system]\label{THM:LIN}
Assume that \ref{ass:1:hpconst}--\ref{ass:3:compat} and \ref{ass:4:data} hold, and let $u \in \UR$ with associated state $(\phi, \mu, w )= \S(u)$
\juerg{be given}.
Then, for every $h \in L^2(Q)$, the linearized system \Lin\ admits a unique solution $(\xi,\eta,\zeta) \in {\cal X} $, 
where $\cal X$ is the Banach space introduced by \eqref{def:X}.
Furthermore, there exists \juerg{some $K_7>0$, which depends only} on 
\gianni{the structure of the system and an upper bound for the norm of $f$ and those of the initial data}, such~that 
\begin{align}
	& \non
	\norma{\xi}_{\H1 \Vp \cap \L\infty V \cap \L2 W}
	+ \norma{\eta}_{\L2 V}
	\\ & \quad 
	+ \norma{\zeta}_{ \H2 \Vp \cap \pier{{}\W{1,\infty} V \cap \H1 W{}}}
%	\last{	+ \norma{\zeta}_{\H2{\last{H}}\cap\W{1,\infty}V \cap \H1W}}
	\leq \gianni{K_7 \norma h_{\L2H}}.
	\label{est:linearized}
\end{align}
\end{theorem}

\gianni{%
\Brem
\label{VariationalLin}
Due to the low regularity level given by the definition \eqref{def:X} of~$\calX$,
the above result must refer to a proper variational formulation of the linearized problem.
For instance, \eqref{lin:sys:1} with the homogeneous Neumann boundary condition for $\eta$ has to be read~as
\Beq
  \< \dt\xi , v >
  + \iO \nabla\eta \cdot \nabla v
  + \gamma \iO \xi v
  = 0
  \quad \hbox{\aet, for every $v\in V$}.
  \non
\Eeq
\Erem
}%

\begin{proof}[Proof of Theorem \ref{THM:LIN}]
As the system is linear, the uniqueness of solutions readily follows once \eqref{est:linearized} has been shown for a special solution. 
Indeed, suppose that there are two solutions $(\xi_1,\eta_1,\zeta_1)$ and $(\xi_2,\eta_2,\zeta_2).$ It is then enough to repeat the procedure used below with $\xi= \xi_1-\xi_2$, $\eta =\eta_1-\eta_2$ and $\zeta=\zeta_1-\zeta_2$ to realize that the same as  \eqref{est:linearized} holds with \pier{the \rhs\ equal to} $0$ so that $(\xi_1,\eta_1,\zeta_1)\equiv (\xi_2,\eta_2,\zeta_2)$, i.e., the uniqueness.

\juerg{Since the proof of existence} is standard, we avoid introducing any approximation scheme and just provide formal estimates. 
The rigorous argument can be straightforwardly reproduced, e.g., on a \juerg{suitable} Faedo--Galerkin scheme. 

\step First estimate

\gianni{We aim at proving that
\begin{align}
	\norma{\xi}_{\L\infty V}
	+ \norma{\eta}_{\L2V}
	+ \norma{\zeta}_{\W{1,\infty}H \cap \H1V}
	\leq C \norma h_{\L2H} \,.
	\label{primastimalin}
\end{align}
We preliminarily observe that 
\begin{align}
  \norma{\dt\xi}_{\Lt2\Vp}
  \leq C \bigl( \norma\xi_{\Lt2H} + \norma\eta_{\Lt2V} \bigr)
  \quad \hbox{for every $t\in(0,T]$}\,,
  \label{an:1}
\end{align}
as one immediately sees by multiplying \eqref{lin:sys:1} by any $v\in\Lt2V$
and integrating over $Q_t$ and by parts.
Moreover, we recall \eqref{wp:est} and \eqref{cor:sep}
and observe that the former yields a uniform $L^\infty$ bound for~$\nabla\phi$
since $\Wx{1,6}\emb\Linfty$.
It \an{then} follows that
\Beq
  \normaV{F''(\phi)\xi}
  \leq C \normaV\xi
  \quad \aet \,.
  \label{pier1}
\Eeq
At this point, we are ready to perform the desired estimate.}
We test \eqref{lin:sys:1} by $\eta+\xi$, \eqref{lin:sys:2} by $-\dt\xi+\eta$, \eqref{lin:sys:3} by $\frac b\lambda \dt\zeta$, 
and add the resulting equalities to infer that \juerg{a.e. in $(0,T)$ it holds}
\begin{align*}
	& \norma{ \eta}^2_V
	+ \frac 12 \, \frac d {dt} \, \norma{\xi}^2_V
	+ \gamma \norma{\xi}^2
	+ \frac b {2\lambda} \, \frac d {dt} \, \norma{\dt \zeta}^2
	+ \frac {\kuno b}{\lambda} \, \norma{\nabla \dt \zeta}^2
	+ \frac {\kdue b}{2\lambda} \, \frac d {dt} \, \norma{\nabla \zeta}^2
	\\ & \quad
	= 
	- \gamma \iO \xi \eta 
	+ \iO F''(\phi) \xi \, (\eta - \dt\xi	)
	-  b \iO \dt \zeta \eta
	+ \frac b \lambda \iO h \, \dt \zeta\,,
\end{align*}
\gianni{thanks to a number of cancellations.
Now, the whole \rhs\ can easily be bounded from above~by
\begin{align*}
	&
	\frac 14 \, \norma{\eta}^2_V
	+ C \bigl( \normaV\xi^2 + \norma{\dt\zeta}^2 + \norma h^2)
	- \iO F''(\phi) \xi \, \dt\xi \,,
\end{align*}
and it is clear that \eqref{primastimalin} follows 
upon integrating in time and invoking Gronwall's lemma
\juerg{provided we can} properly estimate the time integral of the last term.
\an{Using also \eqref{an:1} \pier{and \eqref{pier1}}, we} have~that
\Bsist
  && - \intQt F''(\phi) \xi \, \dt\xi
  \leq C \norma{F''(\phi)\xi}_{\Lt2V} \norma{\dt\xi}_{\Lt2\Vp} 
  \non
  \\
  && \quad 
  \leq C \norma\xi_{\Lt2V} \bigl( \norma\xi_{\Lt2H} + \norma\eta_{\Lt2V} \bigr)
  \leq \frac 14 \, \norma\eta_{\Lt2V}^2
  + C \norma\xi_{\Lt2V}^2\,, 
  \non
\Esist
and this is sufficient to conclude.}

\step Second Estimate

\gianni{We now readily deduce from} \pier{\eqref{an:1} that
%and \eqref{lin:sys:3} and elliptic regularity that
\begin{align}
	\norma{\dt\xi}_{\L2\Vp}
%	+ \|\juerg{\xi}\|_{\L2W}
%	+ \norma{\dtt\zeta}_{\L2\Vp}
	\leq C \norma h_{\L2H}\,.
	\label{pier2}
\end{align}
On the other hand, by comparing the terms in \eqref{lin:sys:2} and taking advantage of 
\eqref{primastimalin} and \eqref{pier1}, well-known elliptic regularity results allow us to infer that}
\begin{align}
% 	\norma{\dt\xi}_{\L2\Vp} +
	\|\juerg{\xi}\|_{\L2W}
%	+ \norma{\dtt\zeta}_{\L2\Vp}
	\leq C \norma h_{\L2H}\,.
	\label{pier3}
\end{align}%

\step Third Estimate

\pier{Now, let us rewrite equation~\eqref{lin:sys:3} \last{in terms of the auxiliary variable $z : = \kuno \dt \zeta + \kdue \zeta + \kuno \lambda \xi$. We obtain}
\begin{align*}
&\frac1{\kuno} \dt z - \Delta z
=  h + \frac{\kdue}{\kuno} \dt \zeta - \kuno \lambda \Delta \xi , 
%\\
%&\quad \hbox{in terms of the auxiliary variable }\, z : = \kuno \dt \zeta + \kdue \zeta + \kuno \lambda \xi ,
\end{align*}
and observe that, in view of \eqref{lin:sys:4}--\eqref{lin:sys:5}, $z$ satisfies Neumann homogeneous boundary conditions and null initial conditions. Then, by known parabolic regularity results, \eqref{primastimalin}\last{,} and \eqref{pier3}\last{,} we easily deduce that 
\begin{align*}
\norma{z}_{\H1H\cap\C0V\cap\L2W} \leq C \Bigl\| h + \frac{\kdue}{\kuno} \dt \zeta - \kuno \lambda \Delta \xi \Bigr\|_{\L2H}\, \leq C \norma h_{\L2H}\,.
\end{align*}
Hence, by recalling the definition of $z$ and the already proved bounds
\eqref{primastimalin}, \eqref{pier2}, \last{and} \eqref{pier3}, we arrive at
\begin{align}
\norma{\zeta}_{\H2{\Vp}\cap\W{1,\infty}V \cap \H1W} \leq C \norma h_{\L2H}\,.
\label{pier4}
\end{align}}%
\pcol{Due to the embeddings $V^* \emb W^* $ and  $W \emb H\equiv H^* \emb W^*$, by interpolation we have that $$\H2{\Vp}\cap \H1W \emb \C1H ,$$ whence \eqref{pier4} entails, in particular, that
\begin{align}
\norma{\zeta}_{\C1 H} \leq C \norma h_{\L2H}.\label{pier4bis}
\end{align}}
This concludes the sketch of the proof.
\end{proof}

We now expect that \juerg{-- provided we select the correct Banach spaces --} the linearized system encapsulates the behavior of the \Frechet\ derivative of the solution operator $\S$.
This is stated rigorously in the next theorem, but \juerg{prior to this}, let us introduce the following Banach space:
\begin{align}
    \non
	{\cal Z} &:= \big(\pier{{}\H1 {W^*} \cap \C0H{}} \cap \L2 W \big) \times \L2 H 
	\\ \label{def:Z}	
	&\ \quad{}\times \big(\pier{ \H2{W^*} \cap \C1 H \cap \H1W }\big).
\end{align}
\begin{theorem}[\Frechet\ differentiability of the solution operator]\label{THM:FRE}
\last{Let \pcol{the set of} assumptions
\ref{ass:1:hpconst}--\ref{ass:3:compat} and \ref{ass:4:data} be fulfilled. Then,}
the control-to-state operator $\S$ is \Frechet\ differentiable at \gianni{any $u\in\UR$} 
as a mapping from $L^2(Q)$ into $\cal Z$.
\gianni{Moreover, for $u\in\UR$, the mapping $D\S(u)\in\calL(L^2(Q),\calZ)$ acts as follows: 
for every $h\in L^2(Q)$, $D\S(u)h$ is the unique solution $(\xi,\eta,\zeta)$ 
to the linearized system \Lin\ associated \juerg{with} $h$}.
\end{theorem}

\begin{proof}[Proof of Theorem \ref{THM:FRE}]
\gianni{We fix $u\in\UR$ and first notice that the map $h\mapsto(\xi,\eta,\zeta)$ of the statement
actually belongs to $\calL(L^2(Q),\calZ)$ as a consequence of \an{\eqref{est:linearized}}.
Then, we} 
proceed with a direct check of the claim by showing that 
\begin{align}\label{fre:formal}
	\frac{\norma{\S(u+h) - \S(u) - (\xi,\eta, \zeta)}_{\cal Z}}{\norma{h}_{L^2(Q)}} \to 0 
	\quad 
	\text{as $\norma{h}_{L^2(Q)} \to 0.$}
\end{align}
This will imply both \juerg{the \Frechet\ differentiability of $\S$ in the sense specified in the statement 
and the validity of the identity} $D\S(u)h = (\xi,\eta, \zeta)$.

\juerg{At this place, we remark that the following argumentation will be formal, because of the low regularity of the
linearized variables (recall Remark 4.3)}. \gianni{
Nevertheless, we adopt it for brevity, in order \juerg{to} avoid any approximation,
like a Faedo--Galerkin scheme based on the eigenfunctions of the Laplace operator with homogeneous Neumann boundary conditions
(in~which case, e.g., the Laplacian of the components of  the discrete solution 
could actually be used as test functions).}

\gianni{Without loss of generality, we \pier{may assume 
that $\norma h_{\LQ2}$ is small enough.}
%namely, such that $u+h\in\UR$.
In particular, we owe to the estimates proved for the solutions to the nonlinear problem
corresponding to both $u$ and $u+h$.}
For convenience, let us set
\begin{align*}
	\psi := \phi^h - \phi - \xi,
	\quad 	
	\sigma := \mu^h - \mu- \eta,
	\quad 
	\omega := w^h - w - \zeta,
\end{align*}
with $(\phi^h, \mu^h, w^h):=\S(u + h)$, $(\phi, \mu, w):=\S(u)$, and where $(\xi, \eta, \zeta)$ is the unique solution to \Lin\ associated with $h$.
Due to the previous results, we already know that $(\psi, \sigma, \omega) \in {\cal X} \emb {\cal Z}$ 
and that, by difference, it yields a \juerg{weak} solution to the system
\begin{alignat}{2}
  & \dt\psi - \Delta\sigma + \gamma \psi
  = 0
  \quad && \text{in $Q$},
  \label{fre:sys:1}
  \\
  & \sigma
  = - \Delta\psi + [F'(\phi^h) - F'(\phi) - F''(\phi)\xi] - b \dt \omega
  \quad && \text{in $Q$},
  \label{fre:sys:2}
  \\
  & \dtt \omega - \Delta(\kuno \dt \omega + \kdue \omega) + \lambda \dt\psi
  = 0
  \quad && \text{in $Q$},
  \label{fre:sys:3}
  \\
  & \dn\psi
  = \dn\sigma
  = \dn(\kuno \dt \omega + \kdue \omega )
  = 0 
  \quad && \text{on $\Sigma$},
  \label{fre:sys:4}
  \\
  & \psi(0) =  \omega(0) =  \dt\omega (0) = 0
  \quad && \text{in $\Omega.$}
  \label{fre:sys:5}
\end{alignat}
\Accorpa\Fre {fre:sys:1} {fre:sys:5}
Besides, with the above notation, \eqref{fre:formal} amounts show that 
\begin{align}\label{fre:proof}
	\norma{(\psi, \sigma, \omega)}_{\cal Z} = o (\norma{h}_{L^2(Q)})
	\quad 
	\text{as $\norma{h}_{L^2(Q)} \to 0$}.
\end{align}
\gianni{Moreover, Theorems~\ref{THM:WP} and~\ref{THM:CD:STRONG} entail that}
%\begin{align}
%	& \non 
%	\norma{\phi^h}_{\W{1,\infty} H  \cap \H1{W} \cap \L\infty {\Wx{2,6}} \cap L^\infty(Q)}
%	+\norma{\mu^h}_{L^\infty(Q)}
%	\\ & \quad 
%	+\norma{w^h}_{\H2 H \cap \W{1,\infty} V}
%	\leq
%	K_5, 
%	\label{fre:reg}
%\end{align}	
\begin{align}
	& 
%	\non 
	\norma{\phi^h}_{\H1 V  \cap \L\infty {\Wx {2,6}} }
	+\norma{\mu^h}_{\L\infty V}
%	\\ & \quad 
	+\norma{w^h}_{\H2 H \cap \W{1,\infty} V}
	\leq
	K_1, 
	\label{fre:reg}
\end{align}	
as well as 
\begin{align}	
	\non
	& \norma{\phi^h - \phi}_{\H1 \Vp \cap \L\infty V\cap \L2 W}
	+ \norma{\mu^h-\mu}_{\L2 V}
	\\ & \quad 
	+ \norma{w^h - w }_{\H2 \Vp \cap \W{1,\infty } H \cap \H1 V}
	\leq K_6 \norma{h}_{\L2 H}.
	\label{fre:cd}
\end{align}
Actually, \gianni{for the logarithmic potential in the two-dimensional setting}, 
we also have a stronger version of \eqref{fre:reg} arising as a consequence of Theorem~\ref{THM:SEP}.

Before entering the details, we recall that Taylor's formula yields that
\begin{align}
	\label{taylor}
	F'(\phi^h) - F'(\phi) - F''(\phi)\xi &= F''(\phi) \psi + R^h \,(\phi^h-\phi)^2,
\end{align}
where the remainder $R^h$ is given by
\begin{align*}
	R^h= \int_0^1 F^{(3)}\big( \phi +s (\phi^h-\phi) \big) (1-s)\,ds\,.
\end{align*}
\gianni{Due to \eqref{cor:sep}, we have~that}
\begin{align}
	\norma{R^h}_{L^\infty(Q)} \leq C.
	\label{stimaresto}
\end{align}

\step First estimate

\gianni{We notice that $\psi$ has zero mean value \pier{as can be easily checked by testing \eqref{fre:sys:1} by $1/|\Omega|$ and using \eqref{fre:sys:5}. Hence, we can}
\an{test} \eqref{fre:sys:1} by $\calN\psi$
and \eqref{fre:sys:2} by $-\psi$.
Moreover, we integrate \eqref{fre:sys:3} in time and test the resulting equation by $\frac b \lambda \dt\omega$.
Finally, we sum~up and add the same term 
$\frac{\kuno b}{2\lambda} \frac d{dt}\norma\omega^2=\frac{\kuno b}{2\lambda}\iO\omega\,\dt\omega$ to both sides.}
We obtain~that
\begin{align*}
	& \frac 12 \, \frac d {dt} \, \norma\psi^2_{*}
	+ \gamma \norma\psi^2_*
	+ \norma{\nabla \psi}^2
	+ \frac b {\lambda} \, \norma{\dt \omega}^2
	+ \frac {\kuno b}{2\lambda} \, \frac d {dt} \, \norma\omega^2_V
	\\ & \quad
	= 
	\iO [F'(\phi^h) - F'(\phi) - F''(\phi)\xi]  \psi
	\pier{{}-{}} \frac { b\kdue}\lambda \iO \pier{\nabla (\an{\boldsymbol{1}*} \omega) \cdot \nabla \dt \omega}
	+\frac {\kuno b} {2\lambda} \iO \omega \,\dt \omega .
\end{align*}
\gianni{Since we aim at applying the Gronwall lemma, we should integrate over $(0,t)$ with respect to time.
However, for brevity, we just estimate the first two terms of the \rhs\ obtained by integration 
(the~last one can be trivially handled by the Young inequality)
and avoid writing the integration variable~$s$ in the integrals over~$(0,t)$}.
The first \gianni{one can be controlled by using} the \Holder\ and Young inequalities, \eqref{fre:cd}, 
the continuous embedding $V\emb\Lx4$, 
\gianni{\eqref{taylor}, \eqref{stimaresto}, and the compactness inequality \eqref{compact} as follows}:
\gianni{%
\begin{align*}
	& \intQt [F'(\phi^h) - F'(\phi) - F''(\phi)\xi]  \psi
	= \intQt [F''(\phi) \psi + R^h \,(\phi^h-\phi)^2] \psi
	\\  &\quad 
	\leq 
	C \iot\norma\psi^2 \, ds
	+ C \iot \norma{\phi^h-\phi}^2_4\norma\psi \, ds
	\leq C \iot \norma\psi^2 \, ds
	+ C \iot \normaV{\phi^h-\phi}^4 \, ds
	\\ & \quad 
	\leq 
	C \iot \norma\psi^2 \, ds
	+ \pier{{}C \hskip1pt T \norma{h}_{L^2(Q)}^4{}}
	\last{{}\leq{} }
	\frac 12 \iot \norma{\nabla \psi}^2 \, ds
	+ C \iot \norma\psi_*^2 \, ds
	+  \pier{{} C \norma{h}_{L^2(Q)}^4{}}.
\end{align*}
As for the second term, \gianni{we integrate by parts both in space and time.
By also accounting for the Young inequality, we find~that }
\begin{align*}
	& \pier{{} -	\frac {b\kdue}\lambda \intQt \nabla(\an{\boldsymbol{1}*}\omega) \cdot \nabla\dt\omega
%	\\
%	& \quad 
	= - \frac {b\kdue}\lambda \iO \nabla(\an{\boldsymbol{1}*}\omega)(t) \cdot \nabla\omega(t)
	+ \frac {b\kdue}\lambda \intQt |\nabla\omega|^2}
	\\
	&\quad 
	\leq \frac {\kuno b}{4\lambda} \iO |\nabla\omega(t)|^2
	+ C \iO \Bigl| \iot \nabla\omega \, ds \Bigr|^2
	+ C \intQt |\nabla\omega|^2
	\leq \frac {\kuno b}{4\lambda} \iO |\nabla\omega(t)|^2
	+ C \intQt |\nabla\omega|^2.
\end{align*}
Thus, we can apply the Gronwall lemma and conclude that} 
\begin{align}
	\norma\psi_{\L\infty \Vp \cap \L2 V}
	+ \norma\omega_{\H1 H \cap \L\infty V}
	\leq C \norma{h}^2_{L^2(Q)}.
\label{pier5}
\end{align}

\step Second estimate

We test \eqref{fre:sys:1} by $\psi$, \eqref{fre:sys:2} by \pier{$\Delta \psi$}, and add the resulting equalities to find that 
\begin{align*}
	\frac 12 \, \frac d {dt} \, \norma\psi^2
	+ \norma{\Delta \psi}^2
	+ \gamma \norma\psi^2
	= \iO [F'(\phi^h) - F'(\phi) - F''(\phi)\xi] \Delta \psi 
	\juerg{\,-\,} b \iO \dt \omega \Delta \psi.
\end{align*}
\gianni{%
As above, we only estimate the \rhs\ of the equality obtained by integrating over~$(0,t)$.
By also accounting for the previous estimate, we \pier{have}~that
\begin{align*}
	& \intQt [F'(\phi^h) - F'(\phi) - F''(\phi)\xi] \Delta \psi 
	\juerg{\,-\,} b \intQt \dt \omega \Delta \psi
	\\ & \quad 
	\leq 
	\intQt |F''(\phi)| \, |\psi| \, |\Delta\psi|
	+ \intQt |R^h| \, |\phi^h-\phi|^2 \, |\Delta\psi|
	+ C \iot \norma{\dt\omega} \, \norma{\Delta\psi} \, ds
	\\ & \quad 
	\leq 
	\frac 12 \iot \norma{\Delta \psi}^2 \, ds
	+ C \iot (\norma\psi^2 + \norma{\dt\omega}^2) \, ds
	+ C \norma{h}_{\LQ2}^4
	\\ & \quad 
	\leq 
	\frac 12 \iot \norma{\Delta \psi}^2 \, ds
	+ C \norma{h}_{\LQ2}^4.
\end{align*}
Thus, \pier{owing also} to the elliptic regularity theory, we conclude~that}
\begin{align}
	\norma\psi_{\pier{\L\infty H\cap{}}\L2 W}
	\leq C \norma{h}^2_{L^2(Q)}.
	\label{pier6}
\end{align}

\step Third estimate

Next, we test \eqref{fre:sys:2} by $\sigma$ and, arguing as above, we obtain that
\begin{align}
	\norma{\sigma}_{\L2 H }
	\leq C \norma{h}^2_{L^2(Q)}.
	\label{pier7} 
\end{align}

\step Fourth estimate

\pier{We can now test \eqref{fre:sys:1} by an arbitrary function $v\in \L2W $ and, in view of 
\eqref{pier6} and \eqref{pier7}, easily infer that
\begin{align*}
	\Bigl| \int_0^T \langle \dt \psi , v \rangle_W  \Bigl|
	&\leq \norma{\sigma}_{\L2 H }\norma{\Delta v}_{\L2 H }
	+ \gamma \norma{\psi}_{\L2 H }\norma{v}_{\L2 H }
	\\
	&\leq C \norma{h}^2_{L^2(Q)}\norma{v}_{\L2 W} \quad \hbox{ for all } \, 	v\in \L2W.
\end{align*}
Hence, $\norma{\dt \psi}_{\L2 {W^*}}$ is uniformly bounded by a quantity proportional to $
\norma{h}^2_{L^2(Q)}$, so that from \eqref{pier6} and an interpolation argument we recover 
that 
\begin{align}
	\norma\psi_{\pier{\H1 {W^*}\cap \C0 H\cap{}}\L2 W}
	\leq C \norma{h}^2_{L^2(Q)}.
	\label{pier8}
\end{align}}%

\step Fifth estimate

\pier{Next, we rewrite equation~\eqref{fre:sys:3} \last{in terms of the auxiliary variable $\tau : = \kuno \dt \omega + \kdue \omega + \kuno \lambda \psi$ to obtain}
\begin{align*}
&\frac1{\kuno} \dt \tau - \Delta \tau
= \frac{\kdue}{\kuno} \dt \omega - \kuno \lambda \Delta \psi \last{.}
%\\
%&\quad \hbox{in terms of the auxiliary variable }\, \tau : = \kuno \dt \omega + \kdue \omega + \kuno \lambda \psi .
\end{align*}
Thanks to \eqref{fre:sys:4}--\eqref{fre:sys:5}, it turns out that $\tau$ satisfies Neumann homogeneous boundary conditions and null initial conditions. Then, by virtue of parabolic regularity results along with \eqref{pier5} and \eqref{pier8}, we have that 
\begin{align*}
\norma{\tau}_{\H1H\cap\C0V\cap\L2W} \leq C \Bigl\| \frac{\kdue}{\kuno} \dt \omega - \kuno \lambda \Delta \psi \Bigr\|_{\L2H}\,\leq C \norma{h}^2_{L^2(Q)}\,.
\end{align*}
Therefore, observing that $\kuno \dt \omega + \kdue \omega = \tau - \kuno \lambda \psi $, it follows that both $\omega$ and $\dt \omega$ satisfy (at least) the same estimate as \eqref{pier8}, which yields 
\begin{align}
	\norma\omega_{\pier{\H2 {W^*}\cap \C1 H\cap{}}\H1 W}
	\leq C \norma{h}^2_{L^2(Q)}.
	\label{pier9}
\end{align}}
This concludes the proof \pier{since the estimates \eqref{pier7}--\eqref{pier9} directly lead to \eqref{fre:proof}.}
\end{proof}

\subsection{Adjoint system and first-order optimality conditions}
\label{SUB:FOC}
As a final step, we now introduce a suitable adjoint system to \Ipbl\ in order to recover a more practical form of  the optimality conditions 
for {\bf (P)}. \pier{Let $u \in \Uad$ be given with its associated state $(\phi,\mu,w)$. In a strong formulation, the adjoint system is expressed by the}  {\it backward-in-time} parabolic system
\begin{alignat}{2}
  & - \dt p - \Delta q + \gamma p + F''(\phi) q - \lambda \dt r = \an{\alpha_1} (\phi - \phi_Q)
  \quad && \text{in $Q$},
  \label{adj:sys:1}
  \\
  & q
  = - \Delta p
  \quad && \text{in $Q$},
  \label{adj:sys:2}
  \\ \non
  & -\dt r - \Delta(\kuno r - \kdue (1 \circledast r)) - b q 
  \\ & \quad 
  =
  \an{\alpha_3} (1 \circledast (w- w_Q))
  + \an{\alpha_4} (w(T) - w_\Omega)
  + \an{\alpha_5} (\dt w - w_Q')
  \quad && \text{in $Q$},
  \label{adj:sys:3}
  \\
  & \dn p
  = \dn q
  = \dn(\kuno r - \kdue (1 \circledast r) )
  = 0 
  \quad && \text{on $\Sigma$},
  \label{adj:sys:4}
  \\ 
%  \non
  & p(T)  =  \an{\alpha_2} (\phi(T)- \phi_\Omega) - \lambda\an{\alpha_6} (\dt w(T)- w_\Omega') ,
  \quad 
%  q(T) =  \an{\alpha_4}(w(T) - w_\Omega),
%  \\ & \quad 
  r(T) = \an{\alpha_6} (\dt w(T)- w_\Omega')
  \quad && \text{in $\Omega$},
  \label{adj:sys:5}
\end{alignat}
\Accorpa\Adj {adj:sys:1} {adj:sys:5}
where the convolution product $\circledast$ has been introduced in \eqref{conv:back}.
\pier{Concerning this product, note in particular that $\dt (1 \circledast r) = -r.$}
Let us \gianni{introduce the following \an{shorthand} for the \rhs\ of \eqref{adj:sys:3}},
\begin{align*}
	f_r := \an{\alpha_3} (1 \circledast (w- w_Q))
  + \an{\alpha_4} (w(T) - w_\Omega)
  + \an{\alpha_5} (\dt w - w_Q')\,.
\end{align*}
\gianni{We also notice that the second term} is independent of time. 
Due to the regularity properties in \eqref{regw} and \ref{ass:7:target}, it holds that 
\begin{align}\label{data:bd}
	\norma{f_r}_{\L2 H} \leq C (\norma{w}_{\H2 H \cap \W{1,\infty} V }+1) \leq C.
\end{align}
\an{
Let us remark that the variable $r$ corresponds to the adjoint of the freezing index $w$.
Besides, equation \eqref{adj:sys:3} is of  first-order in time instead of second-order.
However, it is worth pointing out that  \eqref{adj:sys:3} may be rewritten in the 
time-integrated variable $1 \circledast r$ as it holds that $- \dt r = \dtt (1 \circledast r)$.
%1 ⊛ q, then it turns out that −∂tq = ∂tt(1 ⊛ q), and the system (3.35)–(3.38) looks more
%natural.
}

\begin{theorem}[Well-posedness of the adjoint system]\label{THM:ADJ}
\gianni{Let the assumptions} \ref{ass:1:hpconst}--\ref{ass:3:compat} and \ref{ass:4:data}--\ref{ass:7:target} hold\juerg{, 
and let $u \in \Uad$ with associated state $(\phi,\mu,w)\an{=\S(u)}$ be given}.
Then, the adjoint system \Adj\ admits a unique weak solution $(p,q,r)$ such that 
\begin{align*}
	p  &\in \H1 \Vp \cap \L \infty V \cap \L2 W,
	\\
	q  & \in \L2 V,
	\\
	r & \in \H1 H \cap \L\infty V.
\end{align*}
\end{theorem}

\gianni{%
\Brem
\label{VariationalAdj}
Similarly as in Remark~\ref{VariationalLin},
we should here speak of a proper variational formulation.
For instance, \eqref{adj:sys:1} with the homogeneous Neumann boundary condition for $q$ has to be read~as
\Bsist
  && - \< \dt p , v >
  + \iO \nabla q \cdot \nabla v
  + \iO \bigl( \gamma p + F''(\phi) q - \lambda \dt r \bigr) v
  \non
  \\
  && \quad = \iO \an{\alpha_1}(\phi-\phi_Q) v
  \quad \hbox{\aet, for every $v\in V$}.
  \non
\Esist
\Erem
}%

\begin{proof}[Proof of Theorem \ref{THM:ADJ}]
Again, \gianni{for existence,} we proceed formally \pier{but \last{let us}
underline that the following} computations can however be reproduced
in a rigorous framework. 

\step  First estimate

\an{
We test \eqref{adj:sys:1} by $p+q$, \eqref{adj:sys:2} \an{by $ \dt p + (K_5 +1) q$, where $K_5$ is the positive constant arising from \eqref{cor:sep},} \eqref{adj:sys:3} by $- \frac\lambda b \dt r$ 
and add the resulting identities 
{to each other.} \an{Then, we infer~that}
\begin{align}
	& \non
	- \frac 12 \, \frac d {dt} \, \norma{p}_V^2
	+ (K_5 +1) \norma{q}^2
	+ \norma{\nabla q}^2
	+ \gamma \norma{p}^2
	+ \frac \lambda b \, \norma{\dt r}^2
	\\ & \qquad  \non
	{{}-{}} \frac {\kuno \lambda}{\pier{2}b} \, \frac d {dt} \, \norma{\nabla r}^2
	+ {\frac{\kdue \lambda}{b} \iO \nabla(1\bconv r) \cdot \nabla\dt r}
	\\ & \quad  \non
	=
	- \gamma \iO p q
	- \iO F''(\phi)q (p+q)
	+ \lambda \iO \dt r \, p
	+ \an{\alpha_1} \iO (\phi - \phi_Q) (p+q)
	\\ & \qquad \label{est:adj}
	{+ K_5 \iO \nabla p \cdot \nabla q}
	- \frac  \lambda b \iO f_r \, \dt r.
\end{align}
Now, recalling \eqref{cor:sep}, the second term on the \rhs\ can be bounded from above as
\begin{align*}
	- \iO F''(\phi)q (p+q)
	\leq 
	\norma{F''(\phi)}_\infty \norma{p}\norma{q}
	+ \norma{F''(\phi)}_\infty \norma{q}^2
	\leq 
	\Big(\frac 12 + K_5 \Big) \norma{q}^2
	+ C \norma{p}^2,
\end{align*}
and the first term appearing on the right can be absorbed by the corresponding contribution appearing on the left of \eqref{est:adj}. \pier{By the Young inequality,
we see that the} remaining terms on the \rhs\  are bounded above~by
\begin{align*}
	\pier{\frac \lambda{2b}} \, \norma{\dt r}^2
	+\frac 12 \, \norma{\nabla q}^2
	+\frac 14 \, \norma{ q}^2
	+ C (\norma{p}^2_V+ 1 )\,,
\end{align*}
thanks to \eqref{wp:est} and the \pier{estimate \eqref{data:bd}} of~$f_r$.
Next, we integrate over $(t,T)$, for any $t\in(0,T)$, 
and notice that \ref{ass:7:target} provide uniform bounds for $\normaV{p(T)}^2$ and $\normaV{r(T)}^2$ 
using their explicit form given by \eqref{adj:sys:5}.
\gianni{Moreover, we treat the integral deriving from the last term on the \lhs\ \an{of \eqref{est:adj}} as~follows.
With the notation $Q^t:=\Omega\times(t,T)$, we have~that
\Beq
	\frac{\kdue\lambda}{b} \bintQt \nabla(1\bconv r) \cdot \nabla\dt r
	=  - \frac {\kdue\lambda}{b} \iO \nabla(1\bconv r)(t) \cdot \nabla r(t)
	+ \frac {\kdue\lambda}{b} \bintQt \pier{|\nabla r|^2} \,.
	\non
\Eeq
On the other hand, we also have that
\Bsist
  && \Bigl| -\frac {\kdue\lambda}{b} \iO \nabla(1\bconv r)(t) \cdot \nabla r(t) \Bigr|
  \leq \frac {\kuno\lambda}{\pier{4}b} \iO |\nabla r(t)|^2
  + C \iO |\nabla(1\bconv r)(t)|^2
  \non
  \\
  && \quad  \leq \frac {\kuno\lambda}{\pier{4}b} \iO |\nabla r(t)|^2
  + C \bintQt |\nabla r|^2 \,.
  \non
\Esist
Thus, from \gianni {the (backward) Gronwall} lemma and the obvious \pier{subsequent} inequality
\Beq
  \norma{r(t)} \leq C \, \norma{\dt r}_{\LQ2} + C
  \quad \hbox{for every $t\in[0,T]$}\pier,
  \non
\Eeq
we infer~that}
\begin{align*}
	\norma{p}_{\L\infty V}
	+ \norma{q}_{\L2V}
	+ \norma{r}_{\H1H \cap \L\infty V}
	\leq C.
\end{align*}
}

\step Second estimate

Elliptic regularity theory applied to \eqref{adj:sys:2} then produces
\begin{align*}
	\norma{p}_{\L2 W}
	\leq C.
\end{align*}

\step Third estimate

Finally, it is a standard matter to infer from a comparison argument in \eqref{adj:sys:1}, along with the above estimates, that
\begin{align*}
	\norma{\dt p}_{\L2 \Vp}
	\leq C.
\end{align*}
\gianni{This concludes the (formal) proof of the existence of a solution. 
By performing the same estimates in the case of vanishing \rhs\ and final data,
we see that the solution must vanish, whence uniqueness in the general case \pier{follows} by linearity.}
\end{proof}

Finally, using the adjoint variables, we present the first-order necessary conditions for \pier{an} optimal control $u^*$ solving {\bf (P)}. \pier{In the following, $(p,q,r)$ and $(\xi,\eta,\zeta)$ denote the solutions of the respective adjoint problem and linearized problem, but written in terms of the associated state $(\phi^*,\mu^*,w^*) =\S(u^*)$ that replaces $(\phi,\mu,w)$ in systems \eqref{lin:sys:1}--\eqref{lin:sys:5} and \eqref{adj:sys:1}--\eqref{adj:sys:5}.}

\begin{theorem}[First-order optimality conditions]\label{THM:FOC}
Suppose that \ref{ass:1:hpconst}--\ref{ass:3:compat} and \ref{ass:4:data}--\ref{ass:7:target} hold. 
Let $u^*$ be an optimal control for {\bf (P)} with associated state $(\phi^*,\mu^*,w^*)=\S(u^*)$ and adjoint $(p,q,r)$.
Then, it necessarily \gianni{fulfills} the variational inequality
\begin{align}\label{foc:final}
	 \int_Q (r + \nu u^*)(u-u^*)
		 \geq 0 \quad \gianni{\hbox{ for every $u\in\Uad$}}.
\end{align}
\end{theorem}

\begin{proof}[Proof of Theorem \ref{THM:FOC}]
From standard results of convex analysis, the first-order necessary optimality condition 
for every optimal control $u^*$ of {\bf (P)} is expressed in the abstract \gianni{form~as}
\begin{align*}
	\<D\Jred (u^*), u- u^*>  \geq 0 \quad \forall u \in \Uad,
\end{align*}
\an{where $D\Jred$ denotes the \Frechet\ derivative of the reduced cost functional $\J$.}
As a consequence of the \Frechet\ differentiability of the control-to-state operator established in Theorem \ref{THM:FRE}, 
and the form of the cost functional $\J$ in \eqref{cost}, this entails that any optimal control $u^*$ necessarily fulfills
\begin{align}
	\non
	& {\an{\alpha_1}}\int_Q (\phi^* - \phi_Q)\xi
	+ {\an{\alpha_2}}\iO (\phi^*(T) - \phi_\Omega)\xi(T)
	+ {\an{\alpha_3}}\int_Q (w^* - w_Q)\zeta
	\\ & \quad \non
	+ {\an{\alpha_4}}\iO (w^*(T) - w_\Omega)\zeta(T)
	+  {\an{\alpha_5}} \int_Q (\dt w^* -  w'_Q) \dt \zeta
	\\ & \quad \label{foc:first}
	+  {\an{\alpha_6}} \iO (\dt w^*(T) -  w'_\Omega) \dt \zeta(T)
	+ \nu \int_Q u^* (u-u^*)		
	\geq 0 \quad \forall  u \in \Uad,
\end{align}
where $(\xi,\eta,\zeta)$ is the unique solution to the linearized system 
as obtained from Theorem~\ref{THM:LIN} associated with $(\phi,\mu,w)=(\phi^*,\mu^*,w^*)=\S(u^*)$ and $h=u-u^*$.
Unfortunately, the above formulation is not very useful in numerical applications as it depends on the linearized variables. 
However, with the help of the adjoint variables, playing the role of Lagrangian multipliers, the above variational inequality can be simplified.
In this direction, we test \eqref{lin:sys:1} by $p$, \eqref{lin:sys:2} by $q$, \eqref{lin:sys:3} by $r$, and add the resulting equalities and integrate over time and by parts.
\gianni{More precisely, we should consider the variational formulations of the linearized and adjoint systems
mentioned in Remarks~\ref{VariationalLin} and~\ref{VariationalAdj}
in order to avoid writing some Laplacian that does not exist in the usual sense,
and we should also owe to (well-known) generalized versions of the integration by parts in time.
However, for shortness, we proceed as said above and obtain}
\begin{align*}
	0 = & 
	\int_Q [\dt\xi - \Delta\eta + \gamma \xi] p
%   \\ & \quad  
	+ \int_Q [-\eta - \Delta\xi + F''(\phi)\xi - b \dt \zeta] q
	\\ & \quad  
	+ \int_Q [\dtt \zeta - \Delta(\kuno \dt \zeta + \kdue \zeta) + \lambda \dt\xi - h] r
	\\    
	\separa
	= &
	\int_Q \xi [- \dt p - \Delta q + \gamma p + F''(\phi) q - \lambda \dt r ]
	\\ & \quad  
	+ \int_Q \eta [-\Delta p -q]
	+ \int_Q \dt \zeta[- \dt r - \Delta(\kuno r - \kdue (1 \circledast r)) - b q ]
	\\ & \quad  
	+ \iO [\xi(T) p(T)  + \dt \zeta (T) r(T) + \lambda \xi(T) r(T)]
	- \int_Q h r.
\end{align*}
Using the adjoint system \Adj\ and the associated \gianni{final} conditions, \pier{and integrating by parts as well,} we infer that
\gianni{%
\begin{align*}
  & \intQ r (u-u^*)
  = \intQ hr
  \non
  \\
  \separa
  & = \intQ \xi \, \an{\alpha_1} (\phi^*-\phi_Q)
  + \intQ \dt\zeta \, [
    \an{\alpha_3} (1\bconv(w^*-w_Q)
    + \an{\alpha_4} (w^*(T)-w_\Omega)
    + \an{\alpha_5} (\dt w^*-w'_Q)
  ]
  \non
  \\
  & \quad {}
  + \iO \xi(T) \, [
    \an{\alpha_2} (\phi^*(T)-\phi_\Omega)
    - \lambda \an{\alpha_6} (\dt w^*(T)-w'_\Omega)
  ]
  \non
  \\
  & \quad {}
  + \iO \dt\zeta(T) \, \an{\alpha_6} (\dt w^*(T)-w'_\Omega)
  + \iO \lambda \xi(T) \, \an{\alpha_6} (\dt w^*(T)-w'_\Omega)
  \non
  \\
  \separa
	& = {\an{\alpha_1}}\int_Q (\phi^* - \phi_Q)\xi
	\an{+ {\an{\alpha_2}}\iO (\phi^*(T) - \phi_\Omega)\xi(T)}
	+ {\an{\alpha_3}}\int_Q (w^* - w_Q)\zeta
	\\ & \quad 
	\an{+ {\an{\alpha_4}}\iO (w^*(T) - w_\Omega)\zeta(T)}
	+  {\an{\alpha_5}} \int_Q (\dt w^* -  w'_Q) \dt \zeta
%	+ {\an{\alpha_2}}\iO (\phi^*(T) - \phi_\Omega)\xi(T)
	+  {\an{\alpha_6}} \iO (\dt w^*(T) -  w'_\Omega) \dt \zeta(T),
\end{align*}
}%
so that \eqref{foc:first} entails \eqref{foc:final}, \pier{and this concludes} the proof.
\end{proof}

\begin{corollary}
Suppose the assumptions of Theorem \ref{THM:FOC} are fulfilled, and let $u^*$ be an optimal control with 
associated state $(\phi^*,\mu^*,w^*)=\S(u^*)$ and adjoint $(p,q,r)$.
Then, whenever $\nu >0$, $u^*$ is the $L^2$-orthogonal projection of $- \frac 1 \nu r $ onto $\Uad$.
Besides, we have the pointwise characterization of the optimal control $u^*$ as 
\begin{align}
	\non
	u^*(x,t)=\max \Big\{ 
	u_{\rm min}(x,t), \min\{u_{\rm max}(x,t),-\frac 1{\nu} \, r(x,t)\} 
	\Big\} \quad \hbox{for $a.a. \,(x,t) \in Q.$}
\end{align}
\end{corollary}

%%%%%%%%%%%%%%%%%%%%%%%%%%%%%%%%%%%%%%%%%%%%%%%%%%%%%%%%%%%%%%%%%%%%%%%%%%%%%%

\section*{Acknowledgments}
This research was \pier{partially} supported by the Italian Ministry of Education, 
University and Research (MIUR): Dipartimenti di Eccellenza Program (2018--2022) 
-- Dept.~of Mathematics ``F.~Casorati'', University of Pavia. 
In addition, {PC and AS gratefully acknowledge some other support 
from the MIUR-PRIN Grant 2020F3NCPX ``Mathematics for industry 4.0 (Math4I4)'' and}
their affiliation to the GNAMPA (Gruppo Nazionale per l'Analisi Matematica, 
la Probabilit\`a e le loro Applicazioni) of INdAM (Isti\-tuto 
Nazionale di Alta Matematica). 

%%%%%%%%%%%%%%%%%%%%%%%%%%%%%%%%%%%%%%%%%%%%%%%%%%%%%%%%%%%%%%%%%%%%%%%%%%%%%%

\End{document}

%%%%%%%%%%%%%%%%%%%%%%%%%%%%%%%%%%%%%%%%%%%%%%%%%%%%%%%%%%%%%%%%%%%%%%%%%%%%%%
%%%%%%%%%%%%%%%%%%%%%%%%%%%%%%%%%%%%%%%%%%%%%%%%%%%%%%%%%%%%%%%%%%%%%%%%%%%%%%
\footnotesize
\begin{thebibliography}{99}

\pier{\bibitem{Baio1}		
C. Baiocchi, Sulle equazioni differenziali astratte lineari del primo e del secondo ordine negli spazi di Hilbert, {\it Ann. Mat. Pura Appl. (4)}  {\bf 76} (1967), 233-304.}


\juerg{
\bibitem{BS}
M. Brokate, J. Sprekels,
``Hysteresis and Phase Transitions'', Applied Mathematical Sciences {\bf 121},
Springer, New York, 1996.
}

\juerg{
\bibitem{Cag1}
G. Caginalp,
Stefan and Hele--Shaw type models as asymptotic limits of the phase-field equations,
{\it Phys. Rev. A (3)} {\bf 39} (1989), 5887-5896.
}

\juerg{
\bibitem{Cag2}
G. Caginalp, 
The dynamics of a conserved phase field system: Stefan-like, Hele--Shaw, and Cahn--Hilliard models as asymptotic limits,
{\it IMA J. Appl. Math.} {\bf 44} (1990), 77-94.
}

\juerg{
\bibitem{Cag3}
G. Caginalp, X. Chen,
Convergence of the phase field model to its sharp interface limits,
{\it European J. Appl. Math.} {\bf 9} (1998), 417-445.
}
	
\pcol{\bibitem{CRW}		
C. Cavaterra, E. Rocca, H. Wu, Long-time dynamics and optimal control of a diffuse interface model for tumor growth,
{\it Appl. Math. Optim.}  {\bf 83}  (2021), 739-787.}

\pier{\bibitem{CL}
B. Chen, C. Liu, Optimal distributed control of a Allen--Cahn/Cahn--Hilliard system with temperature,
{\it Appl. Math. Optim.}  {\bf 84}  (2021),  suppl. 2, S1639-S1684.}

\pier{
\bibitem{CGRS}
P. Colli, G. Gilardi, E. Rocca, J. Sprekels, 
Optimal distributed control of a diffuse interface model 
of tumor growth, {\it Nonlinearity} {\bf 30} (2017), 2518-2546.
}

\bibitem{CGRS4}
P. Colli, G. Gilardi, E. Rocca, J. Sprekels,
Well-posedness and optimal control for a Cahn--Hilliard--Oono system
with control in the mass term,  {\it Discrete Contin. Dyn. Syst. Ser. S} {\bf 15} (2022),
2135-2172. 

\an{
\bibitem{CGSS2}
P. Colli, G. Gilardi, A. Signori, J. Sprekels,
Optimal control of a nonconserved phase field model of Caginalp type with thermal memory and double obstacle potential, \pier{{\it Discrete Contin. Dyn. Syst. Ser. S}, Early Access (2023), doi: 10.3934/dcdss.2022210
(see also preprint arXiv:2207.00375 [math.OC], (2022), 1-23).
}}

\bibitem{CGSS3}
P. Colli, G. Gilardi, A. Signori, J. Sprekels,
On a Cahn--Hilliard system with source term and thermal memory,
preprint arXiv:2207.08491 [math.AP], (2022), 1-28.

\pier{
\bibitem{CGS}
P. Colli, G. Gilardi, J. Sprekels, 
Optimal control of a phase field system 
of Caginalp type with fractional operators, 
{\it Pure Appl. Funct. Anal.} {\bf 7} (2022), 1597-1635.
}

\juerg{%
\bibitem{CSS3}
P. Colli, A. Signori, J. Sprekels,
Analysis and optimal control theory for a phase field model of Caginalp type with thermal memory, 
{\it Commun. Optim. Theory} {\bf 4} (2022), doi.org/10.23952/cot.2022.4}

\pcol{%
\bibitem{CSS4}
P. Colli, A. Signori, J. Sprekels,
Optimal control problems with sparsity for tumor growth models involving variational inequalities,
{\it J. Optim. Theory Appl.}  {\bf 194} (2022), 25-58.}

\juerg{
\bibitem{Bio}
E. Dolgi,
What lava lamps and vinaigrette can teach us about cell biology,
{\it Nature} {\bf 555} (2018), 300-302.
}

\pcol{%
\bibitem{GaLaSi}
H. Garcke, K.~F. Lam, A. Signori, 
Sparse optimal control of a phase field tumor model with mechanical effects,
{\it SIAM J. Control Optim.}  {\bf 59}  (2021), 1555-1580.}

\juerg{
\bibitem{GLSS}
H. Garcke, K.~F. Lam, E. Sitka, V. Styles,
A Cahn--Hilliard--Darcy model for tumour growth with chemotaxis and active transport,
{\it Math. Models Methods Appl. Sci.} {\bf 26} (2016), 1095-1148.
}

\pier{
\bibitem{GGM}
A. Giorgini, M. Grasselli, A. Miranville,
The Cahn--Hilliard--Oono equation with singular potential,
{\em Math. Models Methods Appl. Sci.} {\bf 27} (2017), 2485-2510.
}


\juerg{
\bibitem{GN91}
{A.~E. Green, P.~M. Naghdi}, 
{A re-examination of the basic postulates of thermomechanics}, 
{\it Proc. Roy. Soc. London Ser. A} {\bf 432} (1991), 171-194.
}

\juerg{
\bibitem{GN92}
{A.~E. Green, P.~M. Naghdi}, 
{On undamped heat waves in an elastic solid}, 
{\it J. Thermal Stresses} {\bf 15} (1992), 253-264.
}

\juerg{
\bibitem{GN93}
{A.~E. Green, P.~M. Naghdi},
{Thermoelasticity without energy dissipation}, 
{\it J. Elasticity} {\bf 31} (1993), 189-208.%
}

\pier{\bibitem{LioMag}
J.~L. Lions, E. Magenes.
``Non-Homogeneous Boundary Value Problems and Applications'',
Die Grundlehren der Mathematischen Wissenschaft, {\bf 181}, 
Springer-Verlag, Berlin, 1972.}

\pcol{%
\bibitem{MTT}
T.~T.~Medjo, C.~Tone, F.~Tone, 
Maximum principle of optimal control of a Cahn--Hilliard--Navier--Stokes model 
with state constraints,
{\it Optimal Control Appl. Methods}  {\bf 42}  (2021), 807-832.}

\an{%
\bibitem{PS}
V. Pata, M. Squassina, On the strongly damped wave equation, 
{\it Comm. Math. Phys.} {\bf 253} (2005), 511-533.%
}

\pier{\bibitem{Per}		
G. Peralta, Distributed optimal control of the 2D Cahn--Hilliard--Oberbeck--Boussinesq system for nonisothermal viscous two-phase flows,
{\it Appl. Math. Optim.}  {\bf 84}  (2021), suppl. 2, S1219-S1279.}

\juerg{
\bibitem{PG09}
{P.~Podio-Guidugli}, 
{A virtual power format for thermomechanics},
{\it Contin. Mech. Thermodyn.} {\bf 20} (2009), 479-487.
}

\pcol{%
\bibitem{Scar}
L. Scarpa,
Analysis and optimal velocity control of a stochastic convective Cahn--Hilliard equation,
{\it J. Nonlinear Sci.} {\bf 31}  (2021), no.~2, Paper No.~45, 57 pp.}

\pcol{\bibitem{SigTime} 
A. Signori,
Penalisation of long treatment time and optimal control of a tumour growth model 
of {C}ahn--{H}illiard type with singular potential, 
{\it Discrete Contin. Dyn. Syst. Ser. A} {\bf 41} (2021), 2519-2542.%
}%

\bibitem{Simon}
J. Simon,
{Compact sets in the space $L^p(0,T; B)$},
{\it Ann. Mat. Pura Appl.~(4)\/} 
{\bf 146} (1987), 65-96.

\pcol{\bibitem{ST}
J. Sprekels, F. Tr\"oltzsch, 
Sparse optimal control of a phase field system with singular potentials arising in the modeling of tumor growth, 
{\it ESAIM Control Optim. Calc. Var.} {\bf 27} (2021), {suppl., Paper No.~S26, 
27~pp.}}

\pcol{\bibitem{SW}
J. Sprekels, H. Wu,
Optimal distributed control of a Cahn--Hilliard--Darcy system with mass sources,
{\it Appl. Math. Optim.} {\bf 83} (2021), 489-530.}

\end{thebibliography}
